\documentclass[a4paper]{article}

\usepackage[cp1251]{inputenc} 
\usepackage[T2A]{fontenc} 

\usepackage[english]{babel}
\usepackage{amsmath, amsfonts, amssymb, amsthm}
\usepackage{floatflt}
\usepackage{multicol}
\usepackage{amsfonts}
\usepackage{caption}
\usepackage{xcolor}

\usepackage{pdfpages} 
\usepackage[12pt]{extsizes} 

\usepackage{graphicx} 

 \newtheorem{theorem}{Theorem}
\newtheorem{lemma}{Lemma}[section]
\newtheorem{proposition}{Proposition}
\newtheorem{definition}{Definition}

\newtheorem{corollary}{Corollary}[section]

\numberwithin{equation}{section}

\DeclareCaptionLabelSeparator{dot}{. }
\title{Simple closed geodesics on regular tetrahedra in spaces of constant curvature}
\author{Darya Sukhorebska\footnote{
The author is 
funded by the Deutsche Forschungsgemeinschaft (DFG, German Research
Foundation) under Germany's Excellence Strategy EXC 2044-390685587,
Mathematics M\"unster: Dynamics-Geometry-Structure, and
supported by IMU Breakout Graduate Fellowship.}}

\begin{document}
 \date{}
 \maketitle
{\it Abstract}.
   In this survey results on the behavior of 
  simple closed geodesics on regular tetrahedra
 in three-dimensional spaces of constant curvature are presented. \\
 


\tableofcontents

\section{Introduction}

A closed geodesic is called simple if it is not self-intersecting and does not go along itself. 
At the end of the XIX century, working on three body problem, 
Poincar\'e~\cite{Poincare1892}
stated a problem about the existence of geodesic lines on smooth convex two-dimensional surfaces.
 Since then   methods to find closed geodesics 
 on regular surfaces of positive or negative curvature were created.
 In 1927 Birkhoff~\cite{Birkhoff1927} proved that there exists at least one 
 simple closed geodesic on an $n$-dimensional Riemannian manifold homeomorphic to a sphere. 
 In contrast to this,  there are  non-smooth convex closed surfaces in Euclidean space  
 that are free from simple closed geodesics.
From the generalization of the Gauss-Bonnet theorem for polyhedra
 follows a necessary condition for the existence 
of  a simple closed geodesic on   convex polyhedra in $\mathbb{E}^3$. 
 This condition doesn't hold for most convex polyhedra, 
  but it holds for regular polyhedra, in particular regular tetrahedra.

 In this survey we present  results on the behavior of 
  simple closed geodesics on regular tetrahedra
 in three dimensional spaces of constant curvature. 
    D.~Fuchs and E.~Fuchs supplemented and systematized 
  the results on closed geodesics on regular polyhedra in $\mathbb{E}^3$
  (see~\cite{FucFuc07} and~\cite{Fuc09}).
  Protasov~\cite{Pro07} obtained the condition for the existence of simple closed geodesics 
 on an arbitrary tetrahedron in Euclidean space.  
  
Borisenko and Sukhorebska 
  studied simple closed geodesics on  regular tetrahedra
 in three dimensional hyperbolic and spherical spaces 
 (see~\cite{BorSuh2020}, \cite{BorSuh2021} and \cite{Bor2022}). 
In Euclidean space the faces of a tetrahedron have zero
Gaussian curvature, and the curvature of a tetrahedron is concentrated only on its vertices. 
In hyperbolic or spherical space the Gaussian curvature of  faces is $k = -1$ or $1$, and the
curvature of a tetrahedron is determined not only by its vertices, but also by its faces.
In hyperbolic space the planar angle $\alpha$ of a face of a regular tetrahedron 
 satisfies  $0< \alpha < \pi/3$. 
In spherical space the planar angle $\alpha$ satisfies  $\pi/3< \alpha \le 2\pi/3$. 
In both cases the intrinsic geometry of a tetrahedron  
 depends on planar angle.
The behavior of closed geodesics on a regular tetrahedron in a three-dimensional space 
 of constant curvature $k$ differs depending on the sign of $k$. \\

\textit{The author expresses her heartfelt thanks to  Prof. Alexander A. Borisenko 
for setting the problem and  for the valuable discussion.}
 

\section{Historical notes and main results}

 In~\cite{Poincare1892} Henri Poincare studied the properties of the solutions of the three-body problem,
 in particular periodical and asymptotic solutions. 
He found that the key difficulty of this problem can be formulated as an independent problem of
  describing geodesics lines on a convex surfaces.
In~\cite{Poincare1905} Poincare shows existence of a simple (without points of self-intersection)
 closed geodesic on a  convex smooth surface S
that is an embedding of two-dimensional 
sphere into Euclidean space $\mathbb{E}^3$ with induced metric.
He considered the shortest simple
closed curve dividing S into two pieces of equal total Gaussian curvature.
Moreover Poincare stated
a conjecture on the existence of at least three simple 
  closed geodesic on a smooth closed convex two-dimensional surface in $\mathbb{E}^3$.
 Later in 1927, Birkhoff proved
that there exists at least one simple closed geodesic on a $n$-dimensional Riemannian manifold
 homeomorphic to a sphere~\cite{Birkhoff1927}.

In 1929, Lusternik and Schnirelmann~\cite{LustShnir},~\cite{LustShnir47}
 published the proof of Poincare's conjecture. 
However, their proof contained some gaps. 
These were filled in by W. Ballmann in 1978~\cite{Ballmann1978}
 and independently by I. Taimanov in 1992~\cite{Taimanov1992}. 
In 1951-1952, Lusternik and Fet~\cite{LustFet1951},~\cite{Fet1952}
 proved the existence of a closed geodesic on $n$-dimensional 
regular closed manifolds.

Using ideas of Birkhoff it was proved that every Riemannian 
metric on a two dimensional sphere carries infinitely  many geometrically distinct closed
geodesics, cf. Franks~\cite{Franks1992} and Bangert~\cite{Banget1993}. 
The methods of the proof are restricted to surfaces.
  The condition of existence an infinitely many closed geodesics on a compact simply-connected manifold 
of arbitrary dimension is more complicated. 
In 1969 D. Gromoll and W. Meyer~\cite{GromollMeyer1969} shows that there always exist infinitely many
distinct periodic geodesies on an arbitrary compact manifold $M$, provided some
weak topological condition holds: if the sequence
of Betti numbers of the free loop space $LM$ of $M$ is unbounded.
Ziller~\cite{Ziller1977} proved that this condition on free loop space holds
 for symmetric spaces of rank~$>1$.
Rademacher~\cite{Rademacher1989} showed that for a $C^4$-regular metric on 
a compact Riemannian manifold
with finite fundamental group there are infinitely many geometrically distinct closed geodesics. 

In 1898 Hadamard~\cite{Hadamard} showed that on a closed surface of negative curvature
any closed curve, that is not homotopic to zero, could be  deformed into the convex curve
of minimal length within  its free homotopy group.
This minimal curve is unique  and it is  a closed geodesic. 
Then it is interesting  to estimate the number of  closed
geodesics, depending on the length of these geodesics, on a compact manifold of negative curvature.
 Huber~\cite{Huber59},~\cite{Huber61} proved that 
 on a complete closed two-dimensional manifold of constant
  curvature $-1$ the number of closed geodesics of length at most $L$ has the order of growth 
 $e^{L}/L$ as $L~\rightarrow~\infty$. 
 For compact $n$-dimensional manifolds of negative curvature this result was generalized by
  Sinai~\cite{Sinai1966}, Margulis~\cite{Margulis1969}, Gromov~\cite{Gromov1976} and others.

 In Rivin's work~\cite{Rivin}, and later in Mirzakhani's work~\cite{Mirz08},  it's proved that 
 on a complete hyperbolic (constant negative curvature) Riemannian surface 
of genus $g$ and with $n$ cusps
the number of simple closed geodesics of length at most $L$ 
is asymptotic  to (positive) constant times $L^{6g-6+2n}$ as $L~\rightarrow~\infty$. 
You can also see~\cite{Rivin2005},~\cite{ErlandssonSouto2022} for details.  

Theorems about geodesic lines on a convex two-dimensional surfaces play an important role
in geometry ``in the large'' of convex surfaces in spaces of constant curvature. 
Important results on this subject were obtained by Cohn-Vossen~ \cite{Kon-Fosen},
 Alexandrov~ \cite{Alek50}, Pogorelov~\cite{Pogor69}. 
In one of the earliest work   Pogorelov  proved that
 on a closed  convex surface of the Gaussian curvature $\le k,$ $k>0,$
each geodesic of length $< \pi / \sqrt{k}$ realized the shortest path between its endpoint~ \cite{Pog46}.
Toponogov~\cite{Toponog63}
 proved that on  $C^2$-regular closed surface of curvature $\ge k >0$ the length 
of a simple closed geodesic  is at most $2\pi / \sqrt{k}$.
  Vaigant and Matukevich~\cite{VagMatuc} proved that on this surface 
a geodesic of length $\ge 3\pi / \sqrt{k}$ has point of self-intersection.

Geodesics have also been studied on non-smooth surfaces, 
including convex polyhedra in $\mathbb{E}^3$.
Since geodesic is the locally shortest curve then it can not pass through any point for which 
the full angle is less then $2\pi$ (see~\cite{Alek50}).
Gruber~\cite{Gruber1991} showed, 
that in the sense of Baire categories~\cite{ItohRouyerVilcu2015}
most convex surfaces (no regularity required)  do not contain  closed geodesic. 
 Pogorelov~\cite{Pogorelov1949} generalized  Lusternik and Schnirelmann's result 
showing that on any closed convex surfaces 
 there is at least three closed quasi-geodesics. 
Whereas a geodesic has exactly $\pi$ surface angle to either side at each point, 
a quasi-geodesic has at most $\pi$ surface angle to either side at each point. 
Unlike to geodesics, quasi-geodesics can pass through the vertices 
with the full angle $< 2\pi$ on  surface~\cite{Alek78}.

On a convex polyhedron a geodesic has following properties: 
1) it consists of line segments on faces of the polyhedron; 
2) it forms equal angles with edges of adjacent faces; 
3) a geodesic cannot pass through a vertex of a convex polyhedron~\cite{Alek50}.
    Galperin~\cite{Galperin2003} presented a necessary condition for existence 
  a simple closed geodesic on a convex polyhedron in $\mathbb{E}^3$. 
   It is based on a generalization of Gauss-Bonnet theorem for polyhedra.
 The curvature of a convex polyhedron in $\mathbb{E}^3$ is concentrated on its vertices.
   Let $\theta_1,\dots, \theta_n$ be the full angles around the vertices 
   $A_1, \dots, A_n$ of a convex polyhedron.
 The curvature of the vertex $A_i$ is $\omega_i=2\pi-\theta_i$, $i=1,\dots, n$.
 If there is a simple closed geodesic on a convex polyhedron, it is necessary that there is a subset
  $I\subset \{1,2, \dots, n\}$ such that   $$\sum_{i\in I} \omega_i=2\pi.$$
  This condition doesn't hold for most polyhedra, but it holds for  regular polyhedra. 
D. Fuchs and E. Fuchs supplemented  and systematized the results on closed geodesics on
regular polyhedra in three-dimensional Euclidean space (see~\cite{FucFuc07} and~\cite{Fuc09}).
 K. Lawson and others~\cite{Lawson2013} obtain a complete classification of simple closed geodesics on 
 the eight convex  polyhedra (deltahedra) whose faces are all equilateral triangles.

Protasov~\cite{Pro07} obtained a condition for the existence of simple closed geodesics 
 on arbitrary tetrahedron in Euclidean space  and evaluated from above 
 the number of these geodesics
in terms of the difference from $\pi$
 the sum of the angles at a vertex of the tetrahedron.
In particular, it is proved that a simplex has infinitely many  
different simple closed geodesics if and only if all the faces are equal triangles.
A. Akopyan and A. Petrunin~\cite{AkopyanPetrunin2018} showed that
if closed convex surface $M$ in $\mathbb{E}^3$ contains arbitrarily long  simple closed geodesic, 
then $M$ is an  tetrahedron whose faces are equal triangles.

\begin{definition}
A simple closed geodesic on a tetrahedron has  \textit{type $(p, q)$} if
it has $p$~vertices on each of two opposite edges of the tetrahedron, $q$~vertices
on each of other two opposite edges, and $(p + q)$~vertices on each of the remaining
two opposite edges.
\end{definition}

On a regular tetrahedron in Euclidean space, for each ordered pair of
coprime integers $(p, q)$ there exists a whole class of simple closed geodesics of type $(p,q)$, 
up to the isometry of the tetrahedron. 
On the development of the tetrahedron geodesics in each class are parallel to one another.
Furthermore, in each class there is a simple close geodesic passing through the midpoints of two pairs
of opposite edges of the tetrahedron~\cite{BorSuh2020}.

O'Rourke and Vilcu~\cite{RourkeVilcu2022}
 considered simple closed quasi-geodesics on tetrahedra in $\mathbb{E}^3$.

In work~\cite{Davis2017}  Davis and others consider geodesics  which begin and end at vertices
 (and do not touch other vertices) on a regular tetrahedron and cube.
It was proved that a geodesic as above never begins and ends at the same vertex and computed
the probabilities with which a geodesic starting from a given vertex ends at every
other vertex.
Fuchs~\cite{Fuchs2016} obtain similar results for regular octahedron and icosahedron
 (in particular, such a geodesic never ends at the point where it begins).

 Denote by $M^n_k$ a simply-connected complete Riemannian $n$-dimensional manifold
  of constant curvature $k\in \{-1,  0, 1\}$.
  A polyhedron in $M^3_k$ is a surface  obtained by
 gluing finitely many geodesic polygons from $M^2_k$. 
 In particular, a regular tetrahedron in $M^3_k$  is a closed convex polyhedron 
all of whose faces are regular 
 geodesic triangles from $M^2_k$ and all vertices are regular trihedral angles.
 From Alexandrov's gluing theorem~\cite{Alek78} it follows 
 that the polyhedron in $M^3_k$ with the induced metric is a compact  Alexandrov surface $A(k)$  
 with the curvature bounded below by $k$. 
Note, that in $\mathbb{E}^3(M^3_0)$ the curvature of a tetrahedron is concentrated only on its vertices. 
In hyperbolic or spherical space, the Gaussian curvature of  faces is $k = -1$ or $1$ respectively,
 and the curvature of a tetrahedron is determined not only by its vertices, but also by its faces.
  
 In~\cite{RouyerVilcu2015}  Rouyer and Vilcu studied the existence or non-existence
of simple closed geodesics on most 
(in the sense of Baire category~\cite{ItohRouyerVilcu2015}) 
Alexandrov surfaces.
  In particular   it was proved that most surfaces in $A(-1)$
 have infinitely many, pairwise disjoint, simple closed geodesics,
 and most surfaces in $A(1)$ have no simple closed geodesics.

As we say before, on a regular tetrahedron in Euclidean space $\mathbb{E}^3$,
 for each ordered pair of coprime integers $(p, q)$ 
 there exists infernally many  simple closed geodesics of type $(p,q)$, 
 that are parallel each other on the development of the tetrahedron. 
 It's follows from the fact, that 
 the development of a tetrahedron along the geodesic is contained in
the standard triangular tiling of the plane.
Moreover, the vertices of the tiling can be  labeled  in such a way that for any development the
labeling of vertices of the tetrahedron  matches the labeling of vertices of the tiling.
This is something that  holds only for regular tetrahedron and only in $\mathbb{E}^3$~\cite{FucFuc07}.
In general there is no tiling of a plane by regular triangles. 

In spherical space $\mathbb{S}^3$  the planar angle $ \alpha$ of the face 
satisfies  $\pi/3< \alpha \le 2\pi/3$. 
The intrinsic geometry of a tetrahedron depends on $\alpha$.
If the planar angle $\alpha=2\pi/3$, then the tetrahedron  
 is a unit  two-dimensional sphere.
 Hence there are infinitely many simple closed geodesics on it 
 and they are great circles of the sphere. 
 In the following we consider  $\alpha$ such that   $\pi/3 < \alpha < 2\pi/3$.
 In \cite{BorSuh2021} Borisenko and Sukhorebska  proved that
  on a regular tetrahedron in spherical space 
there exists the finite number of simple closed geodesics.
The length of all these geodesics is less than~$2\pi$. 

For any coprime integers $(p, q)$  it was found the numbers $\alpha_1$ and $\alpha_2$,
  depending on  $p$, $q$ and satisfying the inequalities $\pi/3< \alpha_1 < \alpha_2 < 2\pi/3$, such that \\
1) if $\pi/3< \alpha <\alpha_1$, then 
on a regular tetrahedron in spherical space with the planar angle  $\alpha$ 
 there exists unique simple closed geodesic of type $(p,q)$, up to the rigid motion of this tetrahedron,
and  it passes through the midpoints of two pairs of opposite edges of the tetrahedron; \\
2)  if  $\alpha_2 < \alpha < 2\pi/3$, then 
on a regular tetrahedron  with the planar angle   $\alpha$  
there is not simple closed geodesic of type $(p,q)$.

In~\cite{Bor2022} Borisenko proved necessary and sufficient condition of
existence a simple closed geodesic on a regular tetrahedron in $\mathbb{S}^3$. 
We will consider it in details in Section~\ref{inS3}.

 Unlike in $\mathbb{S}^3$, on a regular tetrahedron in hyperbolic space $\mathbb{H}^3$ 
 there are infinitely many simple closed geodesics. 
 Recall that the  planar angle $\alpha$ of a regular tetrahedron  in $\mathbb{H}^3$
 satisfies  $0< \alpha < \pi/3$. 
 In~\cite{BorSuh2020} Borisenko and Sukhorebska
  proved that on a regular tetrahedron in hyperbolic space 
for any coprime integers $(p, q)$, $0 \le p<q$,  
there exists unique simple closed geodesic $\gamma$
 of type $(p,q)$, up to the rigid motion of the tetrahedron. 
$\gamma$ passes through the midpoints of two pairs of opposite edges of the tetrahedron.
These geodesics of type $(p,q)$
 exhaust all simple closed geodesics on a regular tetrahedron in hyperbolic space.
As a part of the proof it was found a constant $d(\alpha)>0$ for $ \alpha \in (0, \pi/3)$ such that 
the  distances from the vertices of the regular tetrahedron to 
a simple closed geodesic is greater then  $d(\alpha)$.
Note, that this property holds only for simple closed geodesics
 on regular tetrahedra in $\mathbb{H}^3$.
In $\mathbb{E}^3$ or $\mathbb{S}^3$ for any $\varepsilon >0$ 
there is a simple closed geodesic $\gamma$ on a regular tetrahedron such that 
the distance from  a vertex of the tetrahedron to $\gamma$ is $<\varepsilon$. 

Furthermore, in~\cite{BorSuh2020} it was proved, that
 the number of simple closed geodesics of length bounded by $L$ is asymptotic to 
$c(\alpha)L^2$, when $L~\rightarrow~\infty$.
If $\alpha~\rightarrow~0$, then $ c(\alpha)~\rightarrow~c_0>0$.
On the other hand, 
when the planar angle $\alpha$ of a regular tetrahedron in hyperbolic space is zero, 
the vertices of the tetrahedron become cusps. 
Then  the tetrahedron becomes a noncompact surface homeomorphic to a sphere
with four cusps, with a complete regular Riemannian metric of constant negative
curvature. The genus of this surface is zero.
 In work of Rivin~\cite{Rivin} it was shown that 
 the number of simple closed geodesics on this surface has order of growth $L^2$.
 
 In~\cite{Bor2022} Borisenko proved, that if 
  planar angles of any tetrahedron in hyperbolic space are at most $\pi/4$, 
then for any pair of coprime integers $(p,q)$ there exists a simple
closed geodesic of type $(p,q)$.
 This situation is differ from  Euclidean space, where 
 there are no simple closed geodesics on a generic tetrahedron~\cite{Galperin2003}.

\section{Closed geodesics on a regular tetrahedron in $\mathbb{E}^3$}
Consider a regular tetrahedron $A_1A_2A_3A_4$ in Euclidean space with the edge of length $1$.

Fix the point of a geodesic on a tetrahedron's edge and roll the
tetrahedron along the plane in such way that the geodesic always touches the plane.
The traces of the faces form the \textit{development} of the tetrahedron on a plane
and the geodesic is a line segment inside the development.

\begin{figure}[h]
\begin{center}
\includegraphics[width=140mm]{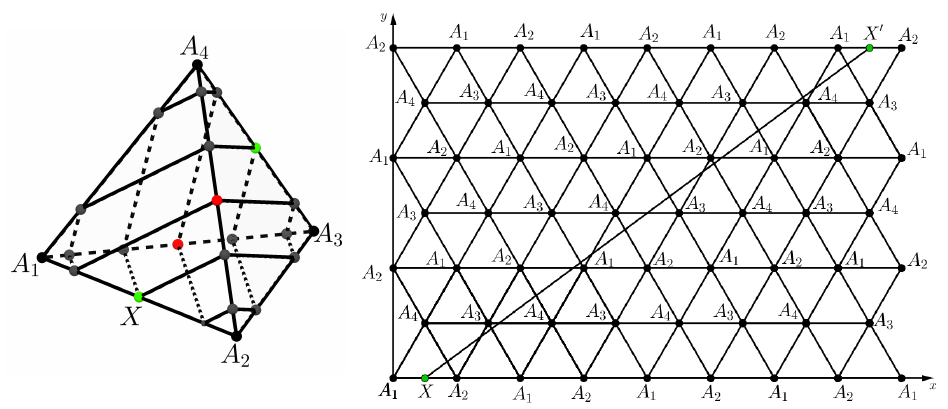}
\caption{ }
\label{evclid_case}
\end{center}
\end{figure}

A development of  a regular  tetrahedron in $\mathbb{E}^3$ is
 a part of the standard triangulation of  Euclidean plane.
Denote the vertices of the triangulation in accordance with the vertices of the tetrahedron 
(see Figure~\ref{evclid_case}).
We introduce a rectangular Cartesian coordinate system with the origin at $A_1$ 
and the $x$-axis along the edge $A_1A_2$ containing $X$.
Then the vertices $A_1$ and $A_2$ has the coordinates $\left( l, k\sqrt{3} \right)$,
and the coordinates of $A_3$ and $A_4$ are $\left( l+1/2, k\sqrt{3} +1/2 \right)$,
where $k, l $ are integers.

Choose two identically oriented edges $A_1A_2$ of the triangulation, 
 which don't belong to the same line.
Take two points $X(\mu, 0)$ and $X'(\mu+q+2p, q\sqrt3)$ on them,
 where $0<\mu<1$ such that the segment $XX'$ doesn't contain any vertex of the triangulation.
The segment $XX'$ corresponds to the simple closed geodesic $\gamma$ of type $(p,q)$
on a regular tetrahedron in Euclidean space.
If $(p,q)$ are coprime integers then $\gamma$ does not repeat itself.
On a tetrahedron $\gamma$ has 
$p$~vertices on each of two opposite edges of the tetrahedron, $q$~vertices
on each of other two opposite edges, 
and $(p + q)$~vertices on each of the remaining two opposite edges, 
so $\gamma$ has type $(p,q)$.

The length of $\gamma$ is equal 
\begin{equation}\label{length_evcl_geod}
L = 2 \sqrt{p^2+pq+q^2}.
\end{equation}

Note, that the segments of a geodesic  lying on the same face of the tetrahedron 
are parallel to each other.
It follows that any closed geodesic on a regular tetrahedron in Euclidean space
does not have points of self-intersection.

If $q=0$ and $p=1$, then geodesic  consists of four segments  
 that consecutively intersect four edges of the tetrahedron,
  and doesn't go through the one pair of  opposite edges.

\begin{theorem}\label{main_th_euclid}
\textnormal{1)}~\textnormal{\cite{FucFuc07} }
 On a regular tetrahedron in Euclidean space, for each ordered pair of
coprime integers $(p, q)$ there exists a whole class of simple closed geodesics of type $(p,q)$, 
up to the isometry of the tetrahedron. 
On the development of the tetrahedron geodesics in each class are parallel  one another.\\
\textnormal{2)}~\textnormal{\cite{BorSuh2020} }
In each class there is a simple close geodesic passing through the midpoints of two pairs
of opposite edges of the tetrahedron.
\end{theorem}

\begin{proof} 
For each pair of coprime integers $(p, q)$ construct  the segment 
connecting points $X(\mu_0, 0)$ and $X'(\mu_0+q+2p, q\sqrt3)$. 
Chose $\mu_0 \in (0, 1)$  such that 
$XX'$ doesn't contain any vertex of the triangulation.
Then $XX'$ corresponds to the simple closed geodesic $\gamma$ of type $(p,q)$
on a regular tetrahedron in Euclidean space.

Consider the segments parallel to $XX'$. 
They are characterized by the equation 
$$y=\frac{q\sqrt3}{q+2p} (x-\mu). $$
We can change $\mu$ until the line touches a vertex of the tiling. 
Then for each pair $(p,q)$ there are $\mu_1, \mu_2 \in(0,1)$ such that 
$\mu_1\le \mu_0 \le \mu_2$ and 
for all $\mu \in (\mu_1, \mu_2)$ the segment joining
  $X(\mu, 0)$ and $X'(\mu+q+2p, q\sqrt3)$ corresponds 
  to the simple closed geodesic of type $(p,q)$
on a regular tetrahedron.
Therefore the part 1) of the theorem is proved. 

To prove 2) consider the lines 
\begin{equation}\label{gamma1}
\gamma_i: y=\frac{q\sqrt3}{q+2p} (x-\mu_i), i=1,2
\end{equation}
pass through the vertices of the tiling.
It means that there exist the integer numbers $c_1$ and $c_2$ such that  the points
$P_1\left( c_1(q+2p)/2q+\mu_1, c_1 \sqrt{3}/2   \right)$
and 
$P_2\left( c_2(q+2p)/2q+\mu_2, c_2 \sqrt{3}/2  \right)$
are the vertices of the tilling and 
 $\gamma_1$ passes through  $P_1$ and $\gamma_2$ passes through $P_2$.

 Consider the closed geodesic $\gamma_0$ parallel
  to $\gamma$ such that the equation of $\gamma_0$ is
\begin{equation}
y=\frac{q\sqrt3}{q+2p} \left( x-\frac{\mu_1+\mu_2}{2} \right).     \notag
\end{equation}
It  passes through the point 
$$P_0 \left( \frac{c_1+c_2}{2}\; \frac{q+2p}{2q} + \frac{\mu_1+\mu_2}{2} ,  \;\;
\frac{c_1+c_2}{2} \;\frac{\sqrt{3}}{2}  \right). $$

Consider three cases: 
1) both of the points $P_1$ and $P_2$ belong to the line $A_1A_2$;
2)  both of the points $P_1$, $P_2$  belong to the line $A_3A_4$;
3) the point  $P_1$  belongs to the line $A_1A_2$ and
the point $P_2$   belongs to the line $A_3A_4$. 
In each of this cases it is easy to show  that $P_0$ is a midpoint of some edge of the tiling.

Then let us proof that if geodesic passes through the midpoint of the one edge, then it passes 
through the midpoints of two pairs of the opposite edges.
Assume that a closed geodesic $\gamma_0$ passes through the midpoint of the edge $A_1A_2$.
Then the equation of $\gamma_0$ is 
\begin{equation}\label{gamma0}
y= \frac{q\sqrt3}{q+2p} \left(x-\frac{1}{2} \right).
\end{equation}

The vertices $A_3$ and $A_4$ belong to the line $y_v=(2k+1) \sqrt{3}/2$,
and their first coordinate is $x_v=l+1/2$ $( k, l \in \mathbb{Z} )$. 
Substituting the coordinates of the points $A_3$ and $A_4$ to  equation (\ref{gamma0}), 
we get 
\begin{equation}\label{condpqkl}
q(2l-2k-1)=2p(2k+1).
\end{equation}

If $q$ is even then there exist $k$ and $l$ satisfying   equation (\ref{condpqkl}).
It follows that $\gamma_0$ passes through the vertex of the tiling.
It contradicts the properties of $\gamma_0$, therefore $q$ is an odd integer.

The points $X_1\left( 1/2, 0 \right)$ and $X'_1 \left(q+2p+1/2, q\sqrt3 \right)$ 
satisfy  equation~(\ref{gamma0}).
These points are the midpoints of the edge $A_1A_2$ on the tetrahedron.
Suppose that the point $X_2$ is the midpoint of $X_1X'_1$.
Then the coordinates of $X_2$ are 
$\left( q/2+p+1/2, q\sqrt3/2 \right)$.
Substituting $q=2k+1$,
we obtain $X_2\left( k+p+1, (k+1/2)\sqrt3 \right)$.
Since the second coordinate of $X_2$ is $(k+1/2)\sqrt3$, where $k$ is an integer, 
then the point $X_2$ belongs to the line, that contains the vertices $A_3$ and $A_4$. 
Since the first coordinate of $X_2$ is an integer,
 it follows that $X_2$ is the midpoint of the edge $A_3A_4$.

Let $Y_1 \left(q/4+p/2+1/2,\; q\sqrt3/4 \right)$ be the midpoint of $X_1X_2$. 
 Substituting $q=~2k+1$, we obtain
 $Y_1 \left((k+p+1)/2+1/4,\; (k/2+1/4)\sqrt3\right)$.
From  the second coordinate we have that $Y_1$ belongs to the line   
passing in the middle of the horizontal lines $y=k\sqrt3/2$ and $y=(k+1)\sqrt3/2$.
Looking at the first coordinate of $Y_1$, which has $1/4$ added, 
we can see that $Y_1$ is the center of
$A_1A_3$, or $A_3A_2$, or $A_2A_4$, or $A_4A_1$.

In the similar way consider the midpoint
 $Y_2 \left(3q/4+3p/2+1/2, \;3q\sqrt3/4\right)$ of   $X_2X'_1$
Then $Y_2$ is the midpoint of the edge that is opposite to the edge with $Y_1$.
\end{proof}

\begin{corollary}\label{sym_develop}
The development of the tetrahedron obtained by unrolling along
a closed geodesic consists of four equal polygons.
Two adjacent polygons can
be transformed into each other by  rotating them through an angle $\pi$ around the midpoint
of their common edge.
\end{corollary}

\begin{proof}
For any  closed geodesic $\gamma$ we get the equivalent closed geodesic $\gamma_0$
 that passes  through the midpoints of two pairs of the opposite edges on the tetrahedron.
Let the points $X_1$, $X_2$  and $Y_1$, $Y_2$ on $\gamma_0$ be  the midpoints 
of the edges $A_1A_2$, $A_4A_3$  and $A_2A_3$, $A_1A_4$ respectively.
\begin{figure}[h!]
\begin{center}
\includegraphics[width=110mm]{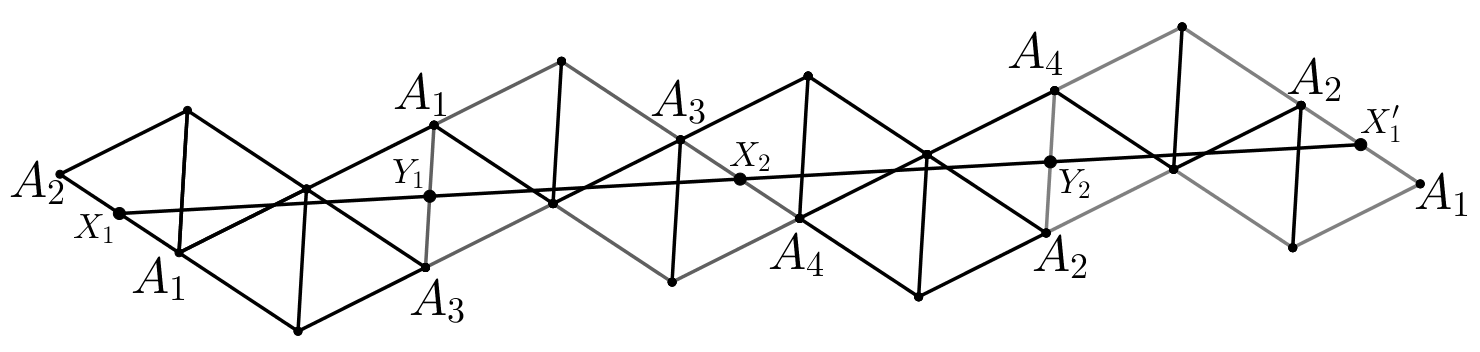}
\caption{  }
\label{geod(2,3)_evcl}
\end{center}
\end{figure}

Consider the rotation of  the regular tetrahedron  through  $\pi$ 
 around the line passing through the points $X_1$ and $X_2$.
This rotation is the isometry of the regular tetrahedron.
The points $Y_1$ and $Y_2$ are swapped.
Furthermore the segment of $\gamma_0$ that starts at $X_1$ on the face  $A_1A_2A_4$  
is mapped to the segment  of $\gamma_0$ that starts from the point $X_1$ on  $A_1A_2A_3$. 
It follows that the  segments  $X_1Y_1$ and  $X_1Y_2$ are swapped.
For the same reason after the rotation the    segments $X_2Y_1$ and $X_2Y_2$ of  $\gamma_0$
are also swapped.

From this rotation we get that the development of the tetrahedron
 along the segment $Y_1X_1Y_2$ of the geodesic is 
central-symmetric  with respect to  the point $X_1$.
And the development along $Y_1X_2Y_2$ is 
central-symmetric  with respect to  $X_2$.

Now consider the rotation of the regular tetrahedron through $\pi$
around the line passing through the points $Y_1$ and $Y_2$.
For the same argument as above we obtain that the development of the tetrahedron 
along the segment $X_1Y_1X_2$ of geodesic is 
central-symmetric  with respect to   $Y_1$,
and the development along the segment $X_2Y_2X_1$ is 
central-symmetric  with respect to  $Y_2$ (see Figure~\ref{geod(2,3)_evcl}).
\end{proof}

\begin{lemma}\label{evcl_dist_vertex_lem}
Let $\gamma$ be a simple closed geodesic of type $(p,q)$ on a regular tetrahedron
in Euclidean space such that
$\gamma$ intersects the midpoints of two pairs of opposite edges.
Then the distance $h$ from the tetrahedron's vertices to $\gamma$ satisfies the inequality 
\begin{equation}\label{evcl_dist_vertex}
h \ge \frac{ \sqrt{3} }{ 4 \sqrt{p^2+pq+q^2} } . 
\end{equation}
\end{lemma}

\begin{proof}
 Suppose  $\gamma$ intersects the edge $A_1A_2$ at the midpoint $X$.
Then  geodesic $\gamma$ is unrolled into the segment $XX'$ lying at the line
  $$y=\frac{q\sqrt3}{q+2p} \left(x-\frac{1}{2} \right).$$
The segment  $XX'$ intersects the edges  $A_1A_2$ at the points 
 $$(x_b,y_b)=\left(  \frac{ 2(q+2p) k + q}{2q},k \sqrt{3} \right),$$ where $k \le q$.
 Since $XX'$ does not pass through   vertices of the tilling,  $x_b$ can not be an integer.
 Hence on the edge  $A_1A_2$ the distance from the vertices 
to the points of $\gamma$ is not less than $1/2q$.
 
Analogically  on the edge  $A_3A_2$ the distance from the vertices of the tetrahedron 
to the points of  $\gamma$ is not less than $1/2p$.

Choose the points $B_1$ at the edge $A_2A_1$ and $B_2$ at the edge $A_2A_3$ such that
 the length $A_2B_1$ is $1/2q$ and  the length $A_2B_2$ equals $1/2p$. 
 Let $A_2H$ be a height  of the triangle $B_1A_2B_2$.
 Then 
 $$ |A_2H| = 
\frac{ \sqrt{3} }{ 4 \sqrt{p^2+pq+q^2} }.$$
The distance $h$ from the vertex $A_2$ to $\gamma$ is not less than  $|A_2H|$.
\end{proof}

The pair of coprime integers $(p,q)$ determines the combinatorical structure of
a simple closed geodesic and 
hence the order of intersections with the edges of the tetrahedron. 

In~\cite{Pro07} the generalization of a simple closed geodesic  on a polyhedron was proposed.
A polyline on a tetrahedron is a curve consisting of line segments which connect
points consecutively on the edges of this tetrahedron.
 \textit{An abstract geodesic}  on a tetrahedron is a closed polyline with the following
properties:\\
1) it does not have points of self-intersection and adjacent segments of it lie on
different faces;\\
2) it crosses more than three edges and does not pass through the vertices of
the tetrahedron.

For any two tetrahedra  we can fix
a one-to-one correspondence between their vertices, and label
the corresponding vertices of the tetrahedra identically.
 Then two closed geodesics
on these tetrahedra are called \textit{equivalent} if they intersect identically labelled edges
in the same order.

 \begin{proposition}\label{allgeod}
\textnormal{\cite{Pro07}}
 For every abstract geodesic $\tilde \gamma$ on a tetrahedron in Euclidean space 
there exists an equivalent simple closed geodesic $\gamma$ on a regular tetrahedron 
in Euclidean space. 
\end{proposition}
 
 A vertex of a geodesic $\gamma$ is called a \textit{link node}
 if it and two neighbouring
vertices of  $\gamma$ lie on the edges of the same vertex $A_i$  of the tetrahedron, and
these three vertices are vertices of the geodesic that are closest to $A_i$.
  
  \begin{proposition}\label{gengeodstr} 
\textnormal{\cite{Pro07}}
Let $\gamma^1_1$ and $\gamma^2_1$ be the   segments of a simple closed geodesic~$\gamma$,
 starting at a link node on a regular tetrahedron,  
let  $\gamma^1_2$ and $\gamma^2_2$ be the next  segments and so on.
 Then for each  $i=2, \dots, 2p+2q-1$ the segments $\gamma^1_i$ and $\gamma^2_i$ 
 lie on the same face of the tetrahedron, and there are no other geodesic points between them. 
 The segments  $\gamma^1_{2p+2q}$ and $\gamma^2_{2p+2q}$ meet at the
second link node of the geodesic.
\end{proposition}

\section{Simple closed geodesics on regular tetrahedra in $\mathbb{S}^3$}\label{inS3}
\subsection{Main definitions and examples.}

A \textit{spherical triangle} is a convex polygon on a unit sphere bounded by three the shortest lines. 
A \textit{regular tetrahedron} $A_1 A_2 A_3 A_4$ in three-dimensional  
spherical space  $\mathbb{S}^3$ is a closed convex polyhedron
such that all its faces are regular spherical triangles and all its vertices are regular trihedral angles.
A planar angle $\alpha$ of a regular tetrahedron in $\mathbb{S}^3$
 satisfies the conditions  $\pi/3 < \alpha \le 2\pi/3$.
Note, than there exist a unique (up to the rigid motion) 
tetrahedron  in spherical space with the given planar  angle.
The   length of the edges is equal to
\begin{equation}\label{a}
a =\text{arccos} \left(  \frac{\cos\alpha}{1-\cos\alpha}  \right), 
\end{equation}
\begin{equation} \label{alim}
\lim\limits_{\alpha\to\pi/3 } a = 0; \;\;\; 
\lim\limits_{\alpha\to\pi/2 } a = \pi/2; \;\;\;
\lim\limits_{\alpha\to2\pi/3 }a = \pi -  cos^{-1}1/3.
\end{equation}

If $\alpha=2\pi/3$, then a tetrahedron is a unit  two-dimensional sphere.
There are infinitely many simple closed geodesics on it.
In the following we assume that  $\alpha$ satisfies   $\pi/3 < \alpha < 2\pi/3$.

Spherical space  $\mathbb{S}^3$  of curvature $1$ is realized as a unite tree-dimensional sphere 
in four-dimensional Euclidean space. 
Hence the regular tetrahedron  $A_1 A_2 A_3 A_4$ is located in an open hemisphere.
Consider  Euclidean space tangent to this hemisphere at the center of  circumscribed sphere 
of the tetrahedron.
A central projection of the hemisphere to this tangent space 
maps the regular tetrahedron from $\mathbb{S}^3$ onto the 
regular tetrahedron in Euclidean tangent space.
A simple closed geodesic  $\gamma$ on  $A_1A_2A_3A_4$ is mapped into  abstract  geodesic
on a regular tetrahedron in $\mathbb{E}^3$.
 Proposition~\ref{allgeod} states  that  there exists 
a simple closed geodesic on a regular tetrahedron 
in  Euclidean space equivalent to this  generalized  geodesic.
It follows, that 
\textit{a simple closed geodesic on a regular tetrahedron in $\mathbb{S}^3$ is
 also characterized  uniquely  by a pair of coprime integers $(p,q)$ 
and has the same combinatorical structure as a closed geodesic on a regular tetrahedron 
 in $\mathbb{E}^3$.}

\begin{lemma}\label{basic_geod}\textnormal{\cite{BorSuh2021}}
1) On a regular tetrahedron with the planar angle $\alpha \in (\pi/3, 2\pi/3)$
 in spherical space there exist three 
different simple closed geodesics of type $(0,1)$.
They coincide under isometries of the tetrahedron. \\
2) Geodesics of type $(0,1)$ exhaust all simple closed geodesics 
on a regular tetrahedron with the planar angle $\alpha \in [\pi/2, 2\pi/3)$
in spherical space.\\
3) On a regular tetrahedron with the planar angle $\alpha \in (\pi/3, \pi/2)$
 in spherical space there exist three 
different simple closed geodesics of type $(1,1)$.
\end{lemma}

\begin{proof}
1) Consider a regular tetrahedron $A_1 A_2 A_3 A_4$ in  $\mathbb{S}^3$
 with planar angle $\alpha \in (\pi/3, 2\pi/3)$. 
 Let $X_1$ and $X_2$ be the midpoints of  $A_1A_4$ and $A_3A_2$,
and $Y_1$, $Y_2$  be the midpoints of  $A_4A_2$ and $A_1A_3$ respectively.
Join these points   consecutively with the segments through the faces. 
Since the points  $X_1$, $Y_1$, $X_2$ and $Y_2$ are midpoints, then 
the triangles  $X_1A_4Y_1$, $Y_1A_2X_2$, $X_2A_3Y_2$ and $Y_2A_1X_1$ are equal.
It  follows that the closed polyline   $X_1Y_1X_2Y_2$  is a simple closed geodesic of type  $(0,1)$
on a regular tetrahedron in spherical space (see Figure~\ref{geod(0,1)_spher}).
Choosing the midpoints of other pairs of opposite edges, 
we can construct other two geodesics of type $(0,1)$
on the tetrahedron.

 \begin{figure}[h!]
\begin{center}
\includegraphics[width=95mm]{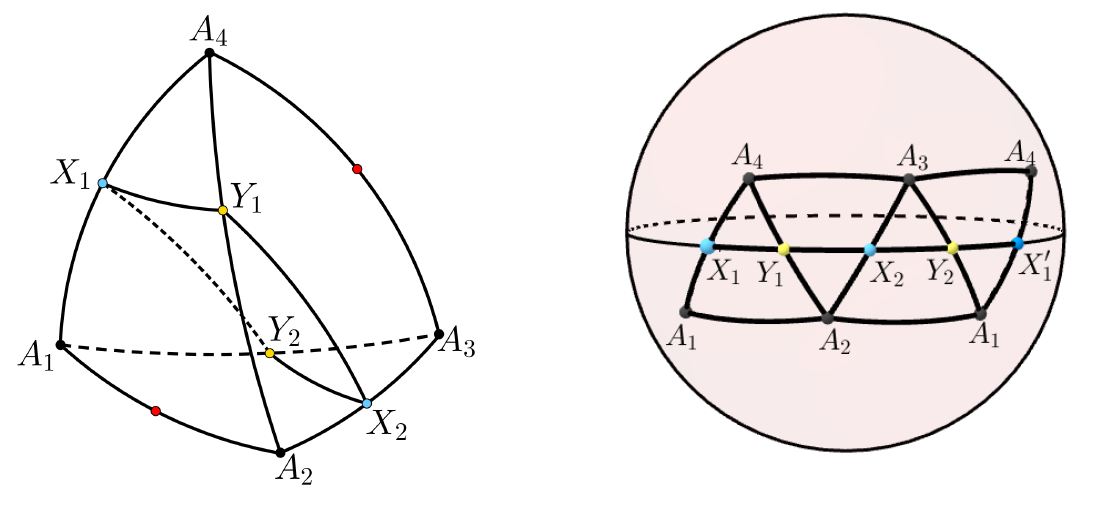}
\caption{ }
\label{geod(0,1)_spher}
\end{center}
\end{figure}

2) Consider a regular  tetrahedron   with planar angle    $\alpha \ge \pi/2$.
Since a geodesic is a line segment inside the development of the tetrahedron, 
then it cannot intersect three edges of the tetrahedron, coming out from the same vertex,  in succession. 

If a simple closed geodesic on the tetrahedron is of type $(p,q)$, 
where  $p~=~q~=~1$ or $1 < p<q$, 
then this geodesic intersect  three  edges, with the common vertex, in succession  (see~\cite{Pro07}). 
Only a simple closed geodesic of  type  $(0,1)$  intersects 
 two tetrahedron's edges, that have  a common vertex,
and doesn't  intersects  the third edge.
It follows that  on a regular tetrahedron in spherical space 
 with planar angle  $ \alpha \in \left[ \pi/2, 2\pi/3 \right) $
 there exist only three simple closed geodesic of type  $(0,1)$ and 
 no other geodesics.

 \begin{figure}[h]
\begin{center}
\includegraphics[width=95mm]{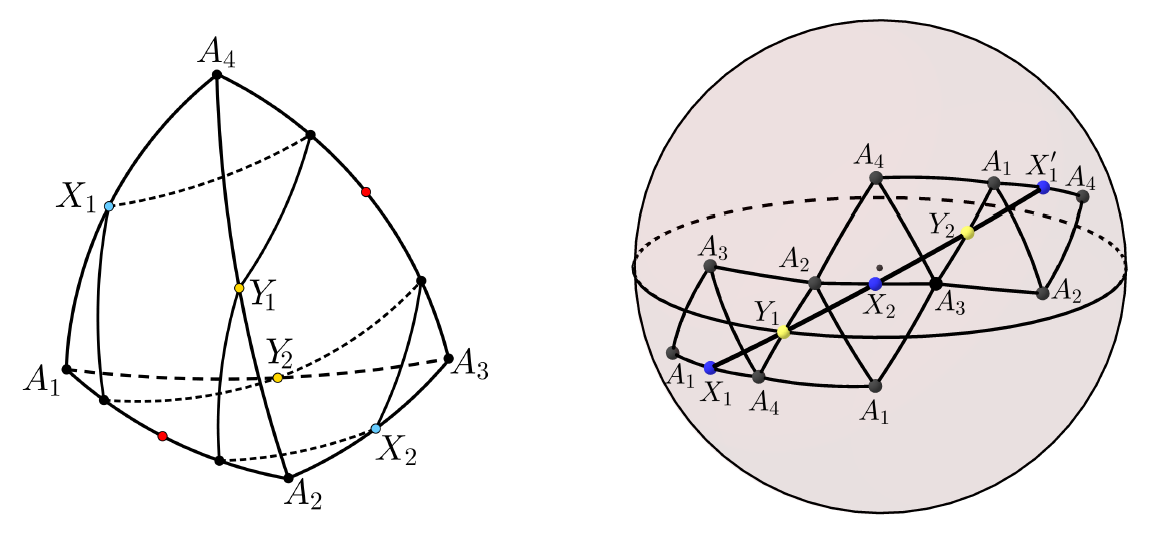}
\caption{ }
\label{geod(1,1)_spher}
\end{center}
\end{figure}

 3) Consider a regular tetrahedron $A_1 A_2 A_3 A_4$ in   $\mathbb{S}^3$
  with planar angle $\alpha\in (\pi/3, \pi/2)$. 
As above, the points $X_1$, $X_2$,  $Y_1$, $Y_2$ are the  midpoints of  
$A_1A_4$, $A_3A_2$, $A_4A_2$ and $A_1A_3$ respectively.

Unfold two adjacent faces  $A_1A_4A_3$ and $A_4A_3A_2$ into the plain and 
draw a geodesic line segment $X_1Y_1$.
Since $\alpha < \pi/2$, then the segment  $X_1Y_1$ is contained inside the development
and intersects the edge  $A_4A_2$  at right angle.
Then unfold another  two adjacent faces $A_4A_1A_2$ and $A_1A_2A_3$ 
and construct the segment $Y_1X_2$.
In the same way join the points  $X_2$ and $Y_2$ within the faces  $A_2A_3A_4$ and  $A_3A_4A_1$,
and join    $Y_2$ and $X_1$ within $A_1A_2A_3$ and $A_4A_1A_2$  (see Figure \ref{geod(1,1)_spher}).
Since  the points  $X_1$, $Y_1$, $X_2$ and $Y_2$ are the midpoints of their edges, 
then the triangles 
$X_1A_4Y_1$, $Y_1A_2X_2$, $X_2A_3Y_2$ and $Y_2A_1X_1$  are equal.
Hence, the segments $X_1Y_1$, $Y_1X_2$, $X_2Y_2$, $Y_2X_1$
 form a simple closed geodesic of type $(1,1)$
on the tetrahedron.

Two other simple closed geodesics of type $(1,1)$
on a tetrahedron can be constructed similarly by connecting 
the midpoints of other pairs of opposite edges of the tetrahedron. 
\end{proof}

In the following we assume that    $\alpha$ satisfying   $\pi/3 < \alpha < \pi/2$.

\subsection{Properties of a simple closed geodesic  on a regular tetrahedron in $\mathbb{S}^3$.}

\begin{lemma}\label{length}
The length of a simple closed geodesic on a regular tetrahedron in spherical space is less than $2\pi$.
\end{lemma}

In~\cite{BorSuh2021} this Lemma was proved using 
Proposition~\ref{gengeodstr}   about the 
construction of a simple closed geodesic on a regular tetrahedron. 
However, Lemma \ref{length} can be considered as the particular case of the result 
proved by A. Borisenko~\cite{Bor2020} 
about the generalization of V. Toponogov theorem~\cite{Toponog63} 
to the case of two-dimensional Alexandrov space.

\begin{lemma}\label{middle_spher}\textnormal{\cite{BorSuh2021}}
On a regular tetrahedron in  spherical space a simple closed geodesic 
intersects midpoints of two pairs of opposite edges.
\end{lemma}
\begin{proof}
Let $\gamma$ be a simple closed geodesic 
on a regular tetrahedron  $A_1 A_2 A_3 A_4$ in   $\mathbb{S}^3$.
As we show above there exists a simple closed geodesic  $\widetilde{ \gamma}$ 
on a regular tetrahedron in Euclidean space such that
 $\widetilde{ \gamma}$  is equivalent to  $\gamma$.
From  Theorem~\ref{main_th_euclid} we  assume  $\widetilde{ \gamma}$ intersects
the midpoints  $\widetilde{X}_1$ and $\widetilde{X}_2$ 
of the edges $A_1A_2$ and $A_3A_4$  on the tetrahedron in $\mathbb{E}^3$.
Denote by  $X_1$ and $X_2$ the vertices of $\gamma$  at the edges $A_1A_2$ and $A_3A_4$
 on the tetrahedron in $\mathbb{S}^3$
such that $X_1$ and $X_2$ are equivalent to the points  $\widetilde{X}_1$ and $\widetilde{X}_2$.

Consider the development of the tetrahedron along $\gamma$ starting from the point  $X_1$
 on a two-dimensional unite sphere.
 The geodesic  $\gamma$ is unrolled into the line segment  $X_1X'_1$ of length less than $2\pi$  inside the development.
 Denote by $T_1$ and $T_2$ the parts of the development along $X_1X_2$ and $X_2X'_1$ respectively.

Let $M_1$ and $M_2$ be   midpoints  
of the edges $A_1A_2$ and $A_3A_4$ respectively 
on the tetrahedron in   $\mathbb{S}^3$.
Rotation by the angle $\pi$ over the line $M_1M_2$ is an isometry of the tetrahedron.
Then the development of the tetrahedron is centrally symmetric with the center $M_2$.

On the other hand, symmetry over $M_2$ swaps the parts $T_1$ and  $T_2$. 
The point  $X'_1$ at the edge $A_1A_2$ of $T_2$ is mapped into the point 
  $\widehat{ X}'_1$  at the edge $A_2A_1$  containing $X_1$ on  $T_1$, and the lengths of
  $A_2X_1$ and $ \widehat{ X}'_1A_1$ are equal. 

The image of the point  $X_1$ on $T_1$ is a point $\widehat{ X}_1$ 
at the edge $A_1A_2$ on $T_2$.
Since $M_2$ is a midpoint of  $A_3A_4$, 
then the symmetry maps the point $X_2$ at  $A_3A_4$ onto
the point $\widehat{X}_2$ at the same edge  $A_3A_4$ 
such that the lengths of $A_4X_2$ and $\widehat{X}_2A_3$ are equal.
Thus, the segment $X_1X'_1$ is mapped into the segment $\widehat{ X}'_1\widehat{ X}_1$ 
inside the development.

\begin{figure}[h]
\centering{\includegraphics[width=80mm]{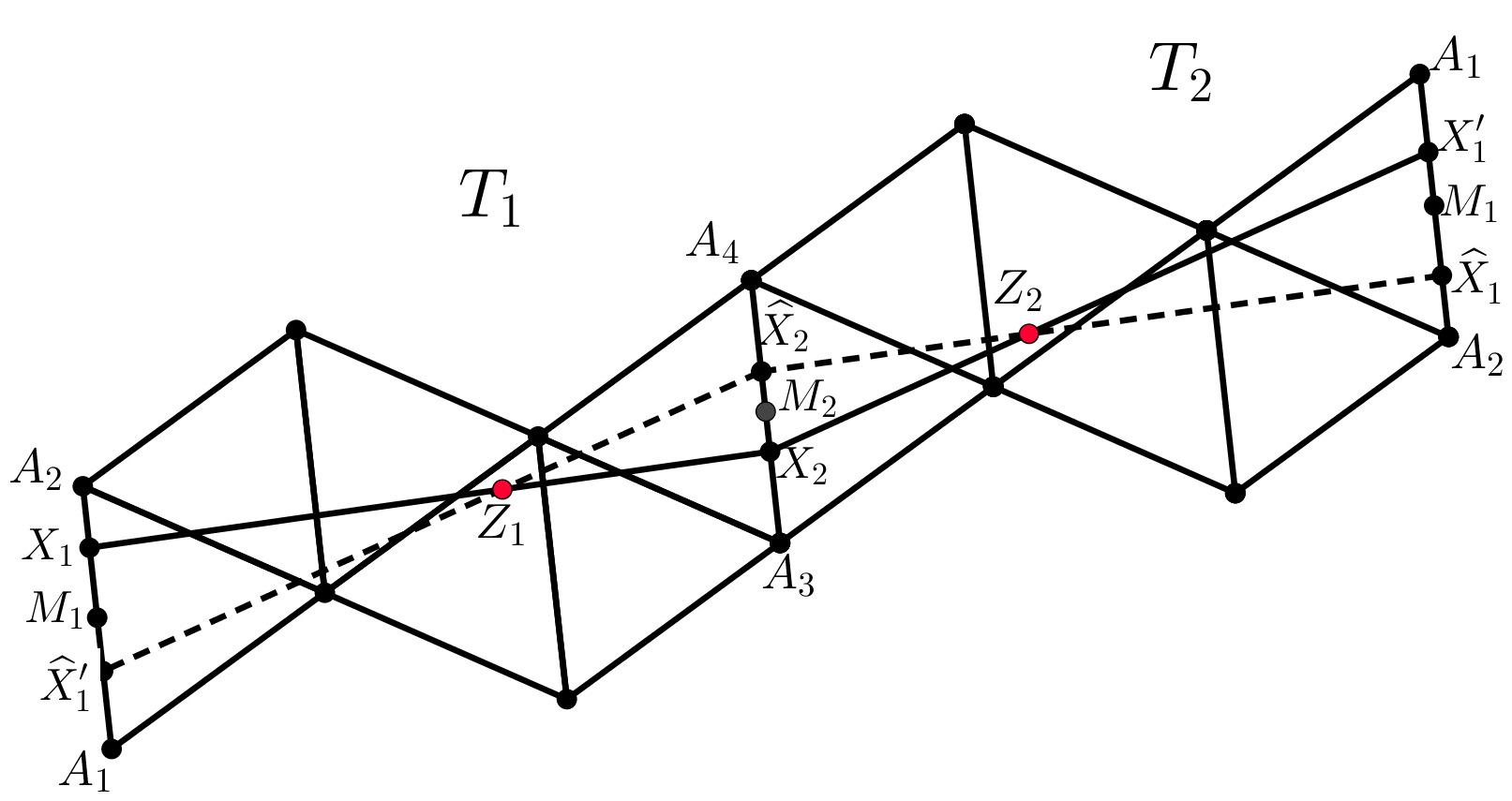} }
\caption{ }
\label{uniqueness_p1}
\end{figure}

Suppose the segments  $\widehat{X}'_1\widehat{X}_2$  and $X_1X_2$ 
intersect at the point $Z_1$ inside $T_1$.
Then the segments $\widehat{ X}_2\widehat{ X}_1$ and $X_2X'_1$ 
intersect at the point $Z_2$ inside $T_2$, 
and the point  $Z_2$  is  centrally symmetric to $Z_1$ with respect to  $M_2$  
(see Figure~\ref{uniqueness_p1}).
Inside the polygon on the sphere 
we obtain two circular arcs $X_1X'_1$ and $\widehat{ X}'_1\widehat{ X}_1$ intersecting in two points.
Therefore  $Z_1$ and $Z_2$ are antipodal points on the sphere and
 the length of the geodesic segment $Z_1X_2Z_2$ equals $\pi$.

Now consider the development of the tetrahedron along $\gamma$ starting from the point  $X_2$. 
This development   also consists of spherical polygons $T_2$ and $T_1$, 
but in this case they are glued by the edge $A_1A_2$ and centrally  symmetric with 
respect to $M_1$.

Similarly  apply the symmetry over $M_1$. 
The segments  $X_2X_1X'_2$ and  $\widehat{ X}_2\widehat{ X}_1\widehat{ X}'_2$
are swapped inside the development.
Since the symmetries over  $M_1$ and over $M_2$ correspond to the same isometry of the tetrahedron,
then  the arcs $X_2X_1X'_2$ and  $\widehat { X}_2\widehat { X}_1\widehat { X}'_2$ also intersect 
at the points  $Z_1$ and $Z_2$.
It follows that the length of geodesic segment  $Z_1X_1Z_2$ is also equal  to $\pi$.
Hence the length of the geodesic $\gamma$ on a regular tetrahedron in spherical space is  $2\pi$,
that  contradicts to   Lemma~\ref{length}.
We get that the segments $\widehat{ X}'_1\widehat{ X}_2$ and $X_1X_2$ on $T_1$
 either don't intersect or coincide.
 
If the $X_1X_2$ and  $\widehat{ X}'_1\widehat{ X}_2$ don't intersect,
 then they form the  quadrilateral $X_1X_2\widehat {X}_2\widehat { X}'_1$  inside $T_1$.
 Since $\gamma$ is closed geodesic, then
  $\angle A_1X_1X_2+\angle A_2\widehat { X}'_1\widehat { X}'_2=~\pi$.
 Furthermore, $\angle X_1X_2A_3+ \angle\widehat {  X}'_1\widehat { X}_2A_4 =\pi$.
 We obtain the convex  quadrilateral on a sphere with  the sum of inner angles  $2\pi$.
It follows that the integral of the Gaussian curvature over the interior of
 $X_1X_2\widehat {X}_2\widehat { X}'_1$
on a sphere is equal zero.
Hence, the segments $X_1X_2$ and   $\widehat { X}'_1\widehat { X}_2$ coincide
 under the symmetry of the development. 
 Then the points $X_1$ and $X_2$ of geodesic $\gamma$  are the midpoints 
  of the edges $A_1A_2$ and $A_3A_4$ respectively.
  
  Similarly it can be proved that  $\gamma$ intersect the midpoints 
  of the second pair of the opposite edges of the tetrahedron. 
  \end{proof}
  
\begin{corollary}\label{uniqueness_spher}\textnormal{\cite{BorSuh2021}}
If two simple closed geodesic on a regular tetrahedron in spherical space intersect
the edges of the tetrahedron in the same order, then they coincide. 
\end{corollary}

 \subsection{An estimation for the angle 
 $\alpha$ for which there is no   simple closed geodesic of type $(p,q)$.} 

  \begin{theorem}\label{necessary_cond}\textnormal{\cite{BorSuh2021}}
  On a regular tetrahedron with the planar angle $\alpha$ in spherical space such that
 \begin{equation}\label{estim_above} 
\alpha > 2\arcsin \sqrt{ \frac{p^2+pq+q^2}{4(p^2+pq+q^2)-\pi^2} }, 
\end{equation}
where $(p,q)$ is a pair of coprime integers, 
there is no simple closed geodesic of type $(p,q)$.
\end{theorem}
 \begin{proof}
Let $A_1 A_2 A_3 A_4$ be a regular tetrahedron  in   $\mathbb{S}^3$  with planar angle  
 $\alpha~\in~(\pi/3, \pi/2)$ and 
 let $\gamma$  be a simple closed geodesic  of type  $(p,q)$ on it.
 
 Each   face  of the tetrahedron is a regular spherical triangle.
 Consider a two-dimensional unit sphere containing the face $A_1A_2A_3$.
Construct   Euclidean plane $\Pi$ passing through the points  $A_1$, $A_2$ and $A_3$.
 The intersection of the sphere with $\Pi$  is a small circle. 
 Draw rays starting at the sphere's center $O$ to the points
  at the spherical triangle $A_1A_2A_3$. 
This defines  the geodesic map between the sphere and the plane  $\Pi$.
 The image of the spherical triangle $A_1A_2A_3$ is 
 the triangle $\widetilde\bigtriangleup A_1A_2A_3$ 
 at the Euclidean plane  $\Pi$.
The edges of $\widetilde \bigtriangleup A_1A_2A_3$  
are the chords joining the vertices of the spherical triangle. 
From (\ref{a}) it follows that the length  $\widetilde{a}$ of an edge of 
 $\widetilde \bigtriangleup A_1A_2A_3$  equals
\begin{equation}\label{tilde_a}
\widetilde {a }= \frac{  \sqrt{ 4\sin^2( \alpha/2 )- 1 }  }{\sin (\alpha/2)   }.
\end{equation}
The segments of the geodesic $\gamma$ lying inside  $A_1A_2A_3$ are mapped into
 the straight line segments inside  $\widetilde \bigtriangleup  A_1A_2A_3$ 
(see Figure \ref{projection_inside}).

 \begin{figure}[h]
\begin{center}
\includegraphics[width=65mm]{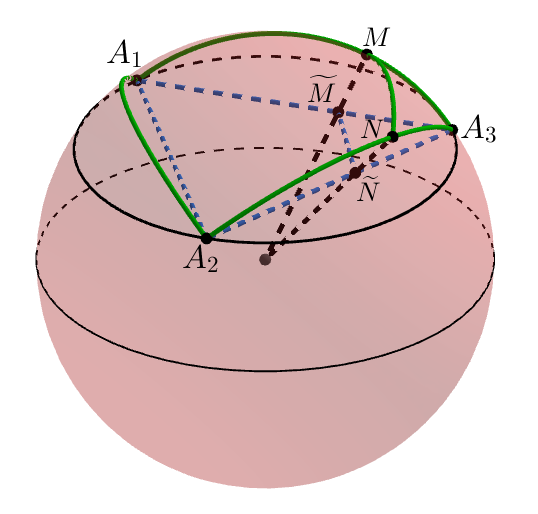}
\caption{ }
\label{projection_inside}
\end{center}
\end{figure}

In the similar way the other tetrahedron faces $A_2A_3A_4$, $A_2A_4A_1$ and $ A_1A_4A_3$ are mapped into 
the plane triangles $\widetilde\bigtriangleup A_2A_3A_4$, 
$\widetilde\bigtriangleup A_2A_4A_1$ and  $\widetilde\bigtriangleup A_1A_4A_3$ respectively. 
 Since the spherical tetrahedron is regular, the constructed plane triangles are equal. 
 We can glue them together identifying the edges with the same labels. 
 Hence we obtain the regular tetrahedron in Euclidean space. 
 Since the segments of $\gamma$ are mapped into
 the straight line segments within the plane triangles, 
 then they form an abstract  geodesic $\widetilde{ \gamma}$  
 on the regular  tetrahedron in $\mathbb{E}^3$,
and  $\widetilde{ \gamma}$ is equivalent to $\gamma$.  

Let us show that the length of $\gamma$ is greater than the length of $\widetilde{ \gamma}$.
 Consider an arc $MN$ of the geodesic $\gamma$ within the face  $A_1A_2A_3$.
  The rays $OM$ and $ON$  intersect the plane $\Pi$ at the points $\widetilde{ M}$ and $\widetilde{ N}$  respectively. 
  The   line segment $\widetilde{ M}$ and $\widetilde{ N}$ lying into
 $ \widetilde\bigtriangleup  A_1A_2A_3$
 is the image of the arc $MN$ under the geodesic map   (see Figure~\ref{projection_inside}).
 Suppose that the length of the arc $MN$ is equal to $2 \varphi $, 
 then the length of the segment  $\widetilde{M}\widetilde{N}$  equals $2 \sin  \varphi $. 
Thus the length of  $\gamma$ on a regular tetrahedron in spherical space is greater than
 the length of its image   $\widetilde{\gamma}$ on a regular tetrahedron in Euclidean space.

From  Proposition~\ref{allgeod} we know that on a regular tetrahedron in  Euclidean space
  there exists a simple closed geodesic $\widehat {\gamma}$ 
    equivalent to  $\widetilde{ \gamma}$. 
On the development of the tetrahedron 
the geodesic  $\widehat {\gamma}$ is a straight line segment, and
the generalized  geodesic   $\widetilde{ \gamma}$ is a polyline, 
then the length of  $\widehat {\gamma}$ is less than the length of $\widetilde{\gamma}$.

This implies that on a regular tetrahedron $A_1 A_2 A_3 A_4$ in $\mathbb{S}^3$
 with planar angle $\alpha$ the length $L_{p,q}$ of 
 a simple closed geodesic   $\gamma$ of type $(p,q)$  is greater than the length of a
 simple closed geodesic $\widehat {\gamma}$ of type $(p,q)$
  on a regular tetrahedron 
 with the edge length  $\widetilde{ a}$ in $\mathbb{E}^3$.
 From the equations  (\ref{length_evcl_geod}) and (\ref{tilde_a})  we get, that
\begin{equation} 
L_{p,q} > 2 \sqrt{p^2+pq+q^2} \frac{\sqrt{4\sin^2(\alpha/2) - 1}}{\sin(\alpha/2)}. \notag
\end{equation}

If $\alpha$ such that the following inequality holds 
\begin{equation}\label{alpha_first} 
2 \sqrt{p^2+pq+q^2}\frac{\sqrt{4\sin^2(\alpha/2) - 1}}{\sin(\alpha/2)} > 2\pi,  
\end{equation} 
then the necessary condition for the existence of a simple closed geodesic of type $(p,q)$
 on a regular tetrahedron with face's angle  $\alpha$  in spherical space is failed. 
Therefore, if 
 \begin{equation}\label{ } 
\alpha > 2\arcsin \sqrt{ \frac{p^2+pq+q^2}{4(p^2+pq+q^2)-\pi^2} },  \notag
\end{equation}
  then there is no simple closed geodesics of type $(p,q)$ 
  on the tetrahedron with planar angle  $\alpha$  in spherical space. 
   \end{proof}

 \begin{corollary}\label{cor_th_1} \textnormal{\cite{BorSuh2021}}
 On a regular tetrahedron in spherical space there exist a finite number of simple closed geodesics.
\end{corollary}

 \begin{proof}
If the integers  $(p,q)$ go to infinity, then 
\begin{equation}
\lim\limits_{p,q \to\infty } 2\arcsin \sqrt{ \frac{p^2+pq+q^2}{4(p^2+pq+q^2)-\pi^2} } 
= 2\arcsin \frac{1}{2} = \frac{\pi}{3}. \notag
\end{equation}
From the inequality  (\ref{estim_above}) we get, that 
for the large numbers  $(p,q)$ a simple closed geodesic of type $(p,q)$
could exist on a regular tetrahedron with the planar  angle $\alpha$ closed to $ \pi/3$ in spherical space.
\end{proof}
 
 The pairs $p=0, q=1$ and $p=1, q=1$ don't satisfy the condition (\ref{estim_above}). 
 Geodesics of this types are described in Lemma~\ref{basic_geod}.

 \subsection{An estimation for the angle 
 $\alpha$  for which there is a simple closed geodesic of type $(p,q)$.} 

In previous sections we assumed that the Gaussian curvature of   faces
of a regular tetrahedron in spherical space  is equal $1$. 
In this case  the length $a$ of the edges of the regular tetrahedron  
was the function of $\alpha$ given by (\ref{a}).
In current section we will assume that 
the faces of the tetrahedron are spherical triangles with the angle $\alpha$
on a sphere of radius $R = 1/a$.
Then the length of the tetrahedron edges equals $1$, and the faces curvature is $a^2$.

Since $\alpha >\pi/3$, then we can write $\alpha =\pi/3 + \varepsilon$, where $\varepsilon > 0$. 
Taking into account Lemma~\ref{basic_geod} we also expect $\varepsilon < \pi/6$.

\begin{theorem}\label{sufficient_cond}\textnormal{\cite{BorSuh2021}}
Let $(p,q)$ be a pair of coprime integers, $0\le p <q$, and let $ \varepsilon$ satisfy
\begin{equation}\label{bar_h_1} 
\varepsilon < \min \left\{
\frac{\sqrt{3}}{4 c_0\sqrt{p^2+q^2+pq}\;
\sum_{i=0}^{\left[ \frac{p+q}{2} \right]+2} \left( c_l (i) + \sum_{j=0}^i c_\alpha (j) \right)};
\frac{1}{8 \cos \frac{\pi}{12}(p+q)^2 }
\right\}, 
\end{equation}
where
\begin{equation}\label{ }
c_0= \frac{ 
3 - \frac{(p+q+2)}{\pi \cos \frac{\pi}{12} (p+q)^2 } - 
16 \sum_{i=0}^{\left[ \frac{p+q}{2} \right]+2} \tan^2 \left( \frac{\pi i}{2 (p+q) } \right) }
{1- \frac{(p+q+2)}{2\pi \cos \frac{\pi}{12} (p+q)^2 } -
8 \sum_{i=0}^{\left[ \frac{p+q}{2} \right]+2} \tan^2 \left( \frac{\pi i}{2 (p+q) } \right) }, \notag
\end{equation} 
\begin{equation}\label{ }
c_l (i) =\frac { \cos \frac{\pi}{12} (p+q)^2 \left( 4+ \pi^2 (2i+1)^2 \right) }{ \left( p+q-i-1 \right)^2 }, \notag
\end{equation}
\begin{equation}\label{ }
c_\alpha (j) =4 \left( 8 \pi (p+q)^2 \cos \frac{\pi}{12} \tan^2 \frac{\pi j}{2 (p+q) } + 1 \right). \notag
\end{equation}
Then on a regular tetrahedron in spherical space with the planar angle $\alpha= \pi/3+\varepsilon$
there exists a unique, up to the rigid motion of the tetrahedron, simple closed geodesic of type $(p,q)$.
\end{theorem}

First let us prove some auxiliary lemmas.

  \begin{lemma}\label{estim_a_above}\textnormal{\cite{BorSuh2021}}
  The edge length of a regular tetrahedron in spherical space of curvature~$1$ satisfies the inequality 
   \begin{equation}\label{a_above_epsilon}
 a <  \pi  \sqrt{2  \cos (\pi/12) } \; \sqrt{  \varepsilon  },        
 \end{equation}
where $\alpha=\pi/3+\varepsilon$ is the planar angle of the face of the tetrahedron.
   \end{lemma}
 
\begin{proof} 
From  (\ref{a}) we have
\begin{equation}\label{sin_a_vs}
\sin a = \frac{ \sqrt{4 \sin^2 (\alpha/2) -1} }{ 2 \sin^2(\alpha/2)  }.  \notag
\end{equation}
Substituting $\alpha =\pi/3+ \varepsilon$, we get
\begin{equation}
\sin a = 
\frac{ \sqrt{\sin(\varepsilon/2) \cos \left(\pi/6 -\varepsilon/2 \right) } }
{ \sin^2 \left( \pi/6 +\varepsilon/2 \right) }. \notag 
\end{equation}
Since $\varepsilon <\pi/6$, then
$$\cos \left(\pi/6 -\varepsilon/2\right) < \cos \pi/12,\;\;\;
\sin \left(  \pi/6 +\varepsilon/2\right) > \sin  \pi/6\;\;\;
\textnormal{and} \;\;\; \sin(\varepsilon/2) < \varepsilon/2.$$
Using this estimations  we obtain
\begin{equation}\label{sin_a}
\sin a < 2 \sqrt{2 \cos (\pi/12) } \;\sqrt{ \varepsilon }. \notag 
\end{equation}
The inequality $a <\pi/2$ implies that $\sin a > (2/\pi) a $. Then
$$a < \pi \sqrt{2 \cos  (\pi/12)  }  \; \sqrt{\varepsilon}. $$
\end{proof}

\begin{figure}[h]
\begin{center}
\includegraphics[width=75mm]{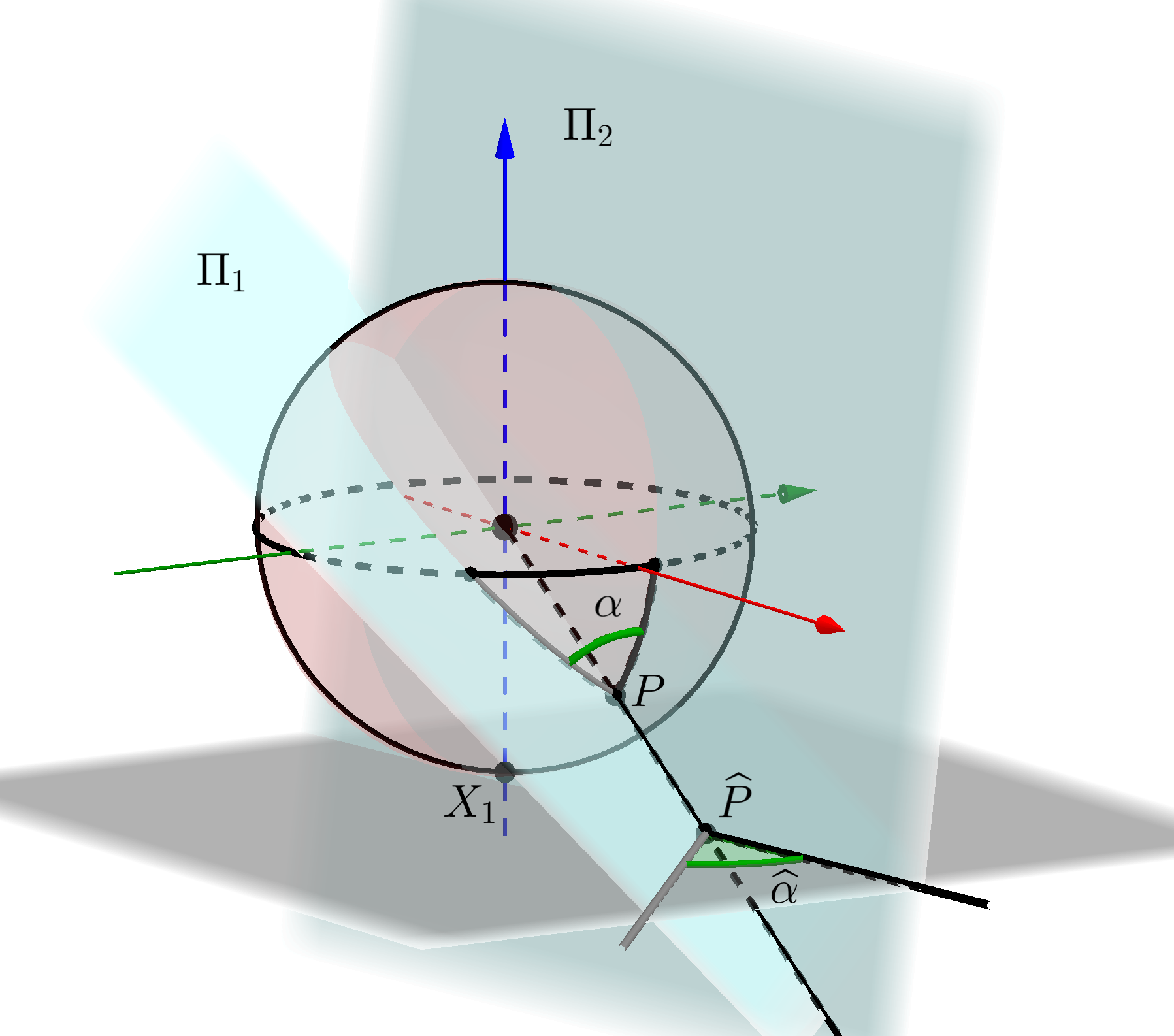}
\caption{ }
\label{angle_projection}
\end{center}
\end{figure}

Consider a parametrization of a two-dimensional sphere $S^2$ of radius $R $ in $\mathbb{E}^3$:
  \begin{equation}\label{spher_th2}
 \begin{cases}
x=R \sin \varphi  \cos \theta \\
y=R  \sin \varphi  \sin \theta \\
z=- R\cos  \varphi 
 \end{cases},  
 \end{equation}
where $ \varphi \in [0, \pi]$, $\theta \in [0, 2\pi)$. 
Let the point $P$ have the coordinates $ \varphi~=~r/R$, $\theta~=~0$, where $ r/R~<~\pi/2$,
and the point $X_1$ correspond to  $ \varphi~=~0$.
Apply a central projection of the hemisphere $ \varphi  \in [0,  \pi/2]$, $\theta \in [0, 2\pi)$  
 onto the tangent plane at  $X_1$ (see Figure~\ref{angle_projection}). 

  \begin{lemma}\label{estim_alpha_above}\textnormal{\cite{BorSuh2021}}
Under the central projection  of the hemisphere of radius $R~=~1/a$ 
onto the tangent  plane at  $X_1$,
 the angle   $\alpha = \pi/3 + \varepsilon$ with the vertex  $P\left(R\sin (r/R), 0, -R\cos (r/R) \right)$ 
on hemisphere is mapped  to the angle    $\widehat{\alpha}_r$ on a plane, which satisfies the inequality 
        \begin{equation}\label{bar_alpha-pi_3}
  \Big| \widehat{\alpha}_r-  \pi/3 \Big| < \pi \tan^2(r/R)+ \varepsilon.   
  \end{equation}
 \end{lemma}

\begin{proof}
Construct the planes $\Pi_1$ and $\Pi_2$  through the center of a hemisphere and 
the point$P\left(R\sin (r/R), 0, -R\cos (r/R) \right)$ :
  \begin{equation}
\Pi_1 : a_1 \cos  (r/R) \: x + \sqrt{1 -a_1^2}\:\: y + a_1 \sin (r/R) \:  z = 0; \notag
  \end{equation}
    \begin{equation}
 \Pi_2 : a_2 \cos  (r/R)\:  x + \sqrt{1-a_2^2} \:\: y + a_2 \sin (r/R)\: z = 0, \notag
  \end{equation}
where 
      \begin{equation}\label{a_1a_2}
     |a_1|, |a_2|  \le 1.
        \end{equation}
If the angle between this two planes   $\Pi_1$ and $\Pi_2$ equals $\alpha$, then
     \begin{equation} \label{cos_alpha_sphere}
\cos \alpha =   a_1a_2 + \sqrt{  (1 -a_1^2)(1 -a_2^2) }.
  \end{equation}

The tangent plane  to $S^2$ at   $X_1$ is given by   $z=-R$.
The planes  $\Pi_1$ and $\Pi_2$  intersect the tangent plane  along the lines, 
that form the angle  $\widehat{\alpha}_r$  (see Figure \ref{angle_projection}),   and
     \begin{equation} \label{cos_alpha_bar}
\cos \widehat{\alpha}_r = \frac{  a_1a_2 \cos^2(r/R)+ \sqrt{ (1 -a_1^2)(1 -a_2^2) }   }
{ \sqrt{1 -a_1^2 \sin^2(r/R)}   \sqrt{1 -a_2^2 \sin^2(r/R) }    }.
  \end{equation}
From the equations (\ref{cos_alpha_sphere}) and (\ref{cos_alpha_bar}) we get 
     \begin{equation}\label{vspomagat_1}
  | \cos \widehat{\alpha}_r- \cos \alpha | < 
  \frac{ |   a_1a_2 \sin^2(r/R) | }{ \sqrt{1 -a_1^2 \sin^2(r/R)  } \sqrt{1 -a_2^2 \sin^2(r/R) }  }. 
  \end{equation}
  Inequalities (\ref{a_1a_2}) and  (\ref{vspomagat_1}) implies that 
\begin{equation}\label{vspomagat_2}
| \cos  \widehat{\alpha}_r - \cos \alpha | < \tan^2 (r/R). 
\end{equation}
It is true that 
\begin{equation}
| \cos \widehat{\alpha}_r- \cos \alpha | = 
\Big|2 \sin  \frac{ \widehat{\alpha}_r - \alpha}{2} \sin \frac{ \widehat{\alpha}_r + \alpha}{2} \Big|  \notag
 \end{equation}
Then  $\alpha > \pi/3$ and $ \widehat{\alpha}_r < \pi $ together with the inequities 
\begin{equation}
 \left |  \sin  \frac{ \widehat{\alpha}_r+ \alpha}{2}  \right| > \sin \frac{\pi}{6}\;\;
 \textnormal{and}\;\;
 \left |  \sin  \frac{ \widehat{\alpha}_r- \alpha}{2} \right|
  > \frac{2}{\pi}  \left | \frac{ \widehat{\alpha}_r- \alpha}{2} \right|,  \notag
  \end{equation}
implies that 
      \begin{equation}\label{bar_alpha_step1} 
     \frac{2}{\pi}  \left | \frac{ \widehat{\alpha}_r - \alpha}{2} \right| < | \cos \widehat{\alpha}_r - \cos \alpha |.  \notag
  \end{equation}
From (\ref{bar_alpha_step1}), (\ref{vspomagat_2}) and   $\alpha = \pi/3 + \varepsilon$ we obtain
\begin{equation}\label{bar_alpha_pi_3_1}
\Big| \widehat{\alpha}_r- \pi/3\Big| < \pi \tan^2(r/R)+ \varepsilon.  \notag
\end{equation}
 \end{proof}

\begin{figure}[h]
\begin{center}
\includegraphics[width=80mm]{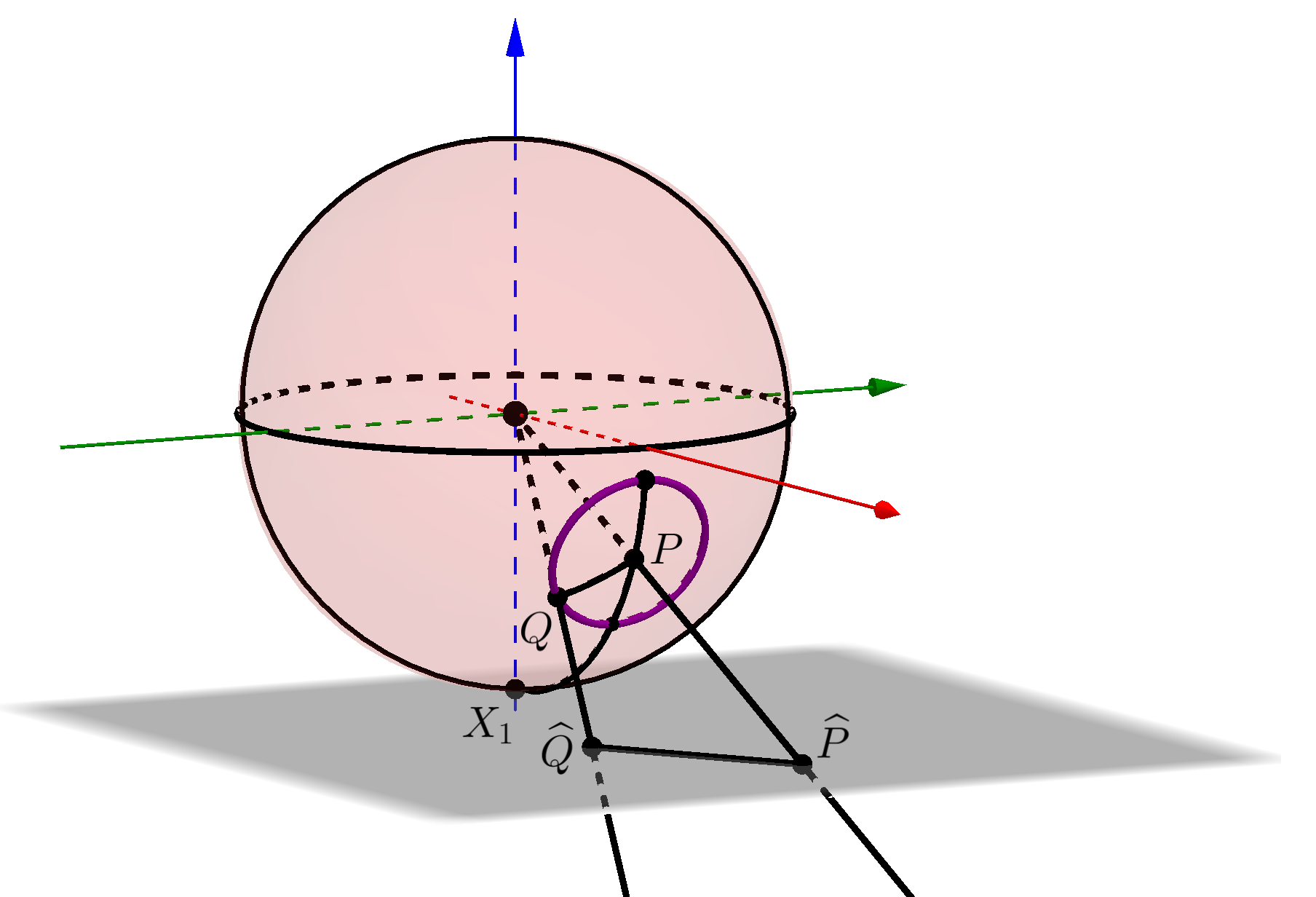}
\caption{ }
\label{length_projection}
\end{center}
\end{figure}

On a sphere  (\ref{spher_th2}) 
let us consider the arc of length one starting at the point  $P$ with the coordinates 
 $ \varphi =r/R, \theta = 0$, where $r/R< \pi/2 $.
 Apply the central projection of this arc to the  plane $z=-R$, 
 which is tangent to the sphere at the point $X_1 ( \varphi  = 0)$ 
 (see Figure~\ref{length_projection}). 

  \begin{lemma}\label{estim_l_above}\textnormal{\cite{BorSuh2021}}
Under the central projection 
 of the hemisphere of radius $R~=~1/a$ onto the tangent  plane at  $X_1$,
  the arc of the length one starting from the point   $P\left(R\sin (r/R), 0, -R\cos (r/R) \right)$ 
is mapped to  the segment of length  $\widehat{l}_r$ satisfying the inequality 
        \begin{equation}\label{bar_l-1}
\widehat{l}_r-1 < \frac {  \cos (\pi/12) \left( 4+  \pi^2  (2r+1)^2  \right) }
{  \left( 1-(2\pi) a (r+1) \right)^2 } \cdot \varepsilon.     
  \end{equation}
 \end{lemma}

\begin{proof}
The point  $P\left(R\sin (r/R), 0, -R\cos (r/R) \right)$   on the sphere  $S^2$ is mapped
 to  $\widehat{P}( R \tan(r/R), 0, -R)$ on the tangent plane $z=-R$.

Take the point $Q(R a_1, R a_2, R a_3)$ on a sphere such that the spherical distance $PQ$ equals~$1$.
Then $\angle~POQ~=~1/R $, where $O $ is a center of the sphere $S^2$ 
(see Figure~\ref{length_projection}).
We obtain  the following conditions for the constants $a_1, a_2, a_3$:
  \begin{equation}\label{cond_1}
  a_1 \sin (r/R) - a_3 \cos  (r/R) = \cos  (1/R);
  \end{equation}
   \begin{equation}\label{cond_2}
  a_1^2+a_2^2+a_3^2 = 1.
  \end{equation}

The central projection into the plane $z=-R$ maps the point $Q $ to the point\\
$\widehat{Q} \left(- \frac{a_1}{a_3} R, - \frac{a_2}{a_3} R, -R \right)$.
The length of $\widehat{P}\widehat{Q}$ equals
\begin{equation}\label{P'Q'}
|\widehat{P}\widehat{Q}| =R \sqrt{ \left( a_1/a_3- \tan (r/R) \right)^2 + a^2_2/a_3^2 }
\end{equation}

Using the  Lagrange multipliers method  to find 
 the local extremum of the length  $\widehat{P}\widehat{Q}$,
we get, that the minimum of $|\widehat{P}\widehat{Q}| $ reaches when 
 $Q$ has the coordinates
  $$ \left( R\sin \left((r-1)/R \right), 0, R\cos \left((r-1)/R \right) \right).$$
  Then
   \begin{equation}\label{ }
 |\widehat{P}\widehat{Q}|_{min} = R     \left|  \tan (r/R)  - \tan \left((r-1)/R \right)\right|  =  
   \frac {R \sin(1/R) }{\cos(r/R)\cos \left((r-1)/R \right)}.   \notag
  \end{equation}
 Note  that $ |\widehat{P}\widehat{Q}|_{min}  > 1$.

  The maximum of $|\widehat{P}\widehat{Q}| $ reaches at the point 
  $Q\left(R\sin\left((r+1)/R \right), 0, R\cos\left((r+1)/R \right) \right)$. 
  This maximum value equals 
   \begin{equation}\label{ }
 |\widehat{P}\widehat{Q}|_{max} = R     \left|   \tan (r/R)   - \tan \left((r+1)/R \right) \right|  =  
   \frac {R \sin(1/R)}{ \cos(r/R) \cos \left((r+1)/R \right)  }.   \notag
  \end{equation}
Since $R =1/a$, then the length $\widehat{l}_r$ of the projection of $PQ$ satisfies 
\begin{equation}\label{bar_l}
\widehat{l}_r< \frac { \sin a }{ a   \cos (a r) \cos \big( a (r+1)\big) }. \notag
\end{equation}
From  $\sin a<a$,  we obtain 
\begin{equation}\label{vspom_bar_l}
\widehat{l}_r-1 < \frac { 2 -  \cos a  - \cos   \big( a(2r+1) \big) }{2   \cos (a r) \cos \big( a (r+1)\big)  }. 
\end{equation}
The equation (\ref{a_above_epsilon}) implies that
\begin{equation}\label{vsp_1}
1 - \cos a = \frac{\sin^2a}{1+ \cos a} \le 8 \cos (\pi/12)\; \varepsilon.
\end{equation}
Similarly from  the inequality  (\ref{a_above_epsilon}) we have
\begin{equation}\label{vsp_2}
1- \cos \left( a (2r+1) \right) \le 2\pi^2 \cos (\pi/12)(2r+1)^2 \; \varepsilon;
\end{equation}   
Estimate the denominator of the   (\ref{vspom_bar_l}) using the inequality
 $\cos x > 1-(2/\pi)x$ where $x <\pi/2$.
Using  (\ref{vsp_1}) and (\ref{vsp_2}), we get
$$
\widehat{l}_r -1 < \frac {   4 \cos (\pi/12) + \pi^2 \cos  (\pi/12) (2r+1)^2}
{    \left( 1- (2/\pi)a \left( r+1 \right) \right)^2  } \cdot \varepsilon.  $$
   \end{proof}


 \begin{proof} \textit{of Theorem~\ref{sufficient_cond}.}
Fix a pair of coprime integers $(p,q)$ such that $0<~p<~q$.
Consider a simple closed geodesic $\widetilde{\gamma}$ of type  $(p,q)$ 
on a regular tetrahedron $\widetilde {A}_1\widetilde{A}_2\widetilde{ A}_3\widetilde{ A}_4$  with 
the  edge of the length~$1$ in $\mathbb{E}^3$.
Assume that $\widetilde{\gamma}$ passes through the midpoints  
$\widetilde{ X}_1$, $\widetilde{ X}_2$ and $\widetilde{ Y}_1$, $\widetilde{ Y}_2$
of the edges $\widetilde{ A}_1\widetilde{A}_2$ and $\widetilde{ A}_3\widetilde{A}_4$
and $\widetilde{A}_1\widetilde{A}_3$, $\widetilde{A}_4\widetilde{ A}_2$ respectively.

Consider the development  $\widetilde{T}_{pq}$ of the tetrahedron 
along  $\widetilde{ \gamma}$ starting from the point $\widetilde{ X}_1$.
The geodesic unfolds to the segment 
$\widetilde{ X}_1\widetilde{ Y }_1\widetilde{ X}_2\widetilde{ Y}_2  \widetilde{X'_1}$ 
inside the development  $\widetilde{T}_{pq}$.
From Corollary~\ref{sym_develop} 
 we know, that the parts of the development along 
geodesic segments $\widetilde{X}_1\widetilde{Y}_1$, $\widetilde{Y}_1\widetilde{X}_2$,
 $\widetilde{X}_2\widetilde{Y}_2$ and $ \widetilde{Y}_2\widetilde{X}'_1$ are equal,
 and any two adjacent polygons can be transformed into each other by a rotation through an angle $\pi$
  around the midpoint of their common edge.

Now consider a two-dimensional sphere $S^2$ of radius $R = 1/a$,
 where  $a$ depends on $\alpha$ according to (\ref{a}).
On this sphere we take the several copies of regular spherical triangles
 with the angle $\alpha \in (\pi/3, \pi/2)$ at vertices.
Fold this triangles up in the same order as  the faces of the Euclidean tetrahedron were unfolded 
along $\widetilde{\gamma}$ into the plane.
In other words, we construct a polygon $T_{pq}$ on a sphere $S^2$
formed by the same sequence of regular triangles
as the polygon $\widetilde{T}_{pq}$ in  $\mathbb{E}^3$.
 Denote the vertices of $T_{pq}$ in accordance with to the vertices of $\widetilde{T}_{pq}$.
By the construction  the spherical polygon $T_{pq}$
 has the same properties of the central symmetry 
as the Euclidean  $\widetilde{T}_{pq}$.
Since the groups of isometries of   regular tetrahedra in $\mathbb{S}^3$
 and in  $\mathbb{E}^3$ are equal, 
then $T_{pq}$ corresponds to the development of a regular tetrahedron 
with the planar angle   $\alpha$ in spherical space.

Denote by $X_1$, $X'_1$ and $X_2$, $Y_1$, $Y_2$ the midpoints of the edges 
$A_1A_2$, $A_3A_4$,  $A_1A_3$,  $A_4A_2$ on $T_{pq}$ respectively. 
These midpoints correspond to the points 
 $\widetilde{X}_1$, $\widetilde{X}'_1$ and $\widetilde{X}_2$, $\widetilde{Y}_1$, $\widetilde{Y}_2$ 
 on the Euclidean development  $\widetilde{T}_{pq}$.
Construct the great circle arcs $ X_1Y_1$, $Y_1X_2$, $ X_2Y_2$  and $ Y_2X'_1$.
The central symmetry of $T_{pq}$ implies that these arcs form the one great arc   $X_1X'_1$
on  $S^2$.
If $\alpha$ such that $X_1X'_1$ lies inside $T_{pq}$, then 
 $X_1X'_1$  correspond  to the  simple closed geodesic of type  $(p,q)$ 
on a regular tetrahedron with the planar angle  $\alpha$ in  $\mathbb{S}^3$.

In what  follows we  consider the part of the polygon  $T_{pq}$ only along $X_1Y_1$ but
we also denote it as $T_{pq}$ for convenience.
This part consists of $p+q$ regular spherical triangles with the edges of length~$1$. 
The  polygon $T_{pq}$ is contained inside the open hemisphere if  
 \begin{equation}\label{a_p+q}
a (p+q) < \pi/2,
\end{equation}
Since $\alpha~=~\pi/3~+~\varepsilon$, then 
 the condition  $(\ref{a_above_epsilon})$ implies that  (\ref{a_p+q}) holds if
\begin{equation}\label{epsilon_estim_1}
\varepsilon < \frac{1}{8 \cos (\pi/12)(p+q)^2 }. 
\end{equation}
In this case the length of the arc $X_1Y_1$ is less than  $\pi/2 a$, 
so   $X_1Y_1$ satisfies the necessary condition from Lemma~\ref{length}.

Apply a central projection of the $T_{pq}$ into the tangent plane  $T_{X_1}S^2$
 at the point $X_1$ to the sphere $S^2$. 
The image of the spherical polygon  $T_{pq}$ on  $T_{X_1}S^2$
is a polygon  $\widehat{T}_{pq}$. 

Denote by  $\widehat{A}_i$ the vertex of $\widehat{T}_{pq}$, 
which is an image of the vertex $A_i$ on $T_{pq}$. 
The arc $X_1Y_1$ maps into the line segment $\widehat{X}_1\widehat{Y}_1$  on   $T_{X_1}S^2$, that
joins  the midpoints of the edges  $\widehat{ A}_1\widehat{ A}_2$ and $\widehat{ A}_1\widehat{ A}_3$.
If for some $\alpha$ the segment $\widehat{ X}_1\widehat{ Y }_1$
 lies inside the polygon  $\widehat{ T}_{pq}$, then 
the arc $X_1Y_1$ is also containing inside  $T_{pq}$ on the  sphere. 

The vector $\widehat{ X}_1\widehat{Y}_1$ equals
     \begin{equation}\label{X_1_overline_Y_1 } 
\widehat{ X}_1\widehat{Y}_1= \widehat{a_0}+\widehat{a_1}+\dots+\widehat{ a_s}+\widehat{a}_{s+1},
 \end{equation}
 where $\widehat{a_i}$ are the sequential vectors of the $\widehat{T}_{pq}$  boundary,  
 $\widehat{a_0}=\widehat{ X_1}\widehat{ A_2}$,  $\widehat{a}_{s+1} =\widehat{ A_1}\widehat{Y_1}$,
and $s=\left[ \frac{p+q}{2} \right]+1$
  (if we take the boundary of $\widehat{ T}_{pq}$ from the other side of
 $\widehat{ X}_1\widehat{ Y}_1$, then $s=\left[ \frac{p+q}{2} \right]$ )
 (see Figure \ref{dev_compare}).
 
On the other hand  at Euclidean plane $T_{X_1}S^2$ 
there exists a development $\widetilde{T}_{pq}$ of a regular Euclidean tetrahedron 
  $\widetilde{A}_1\widetilde{  A}_2\widetilde{ A}_3\widetilde{  A}_4$ with the edge of length~$1$
  along a simple closed geodesic  $\widetilde{\gamma}$.
The development  $\widetilde{T}_{pq}$  is equivalent to $T_{pq}$, 
and then  it's equivalent to  $\widehat{T}_{pq}$.
  The segment  $ \widetilde{  X}_1\widetilde{ Y}_1$ lies inside  $\widetilde{  T}_{pq}$
  and corresponds to the segment of $\widetilde{\gamma}$   (see Figure~\ref{dev_compare}). 

 \begin{figure}[h]
\begin{center}
\includegraphics[width=65mm]{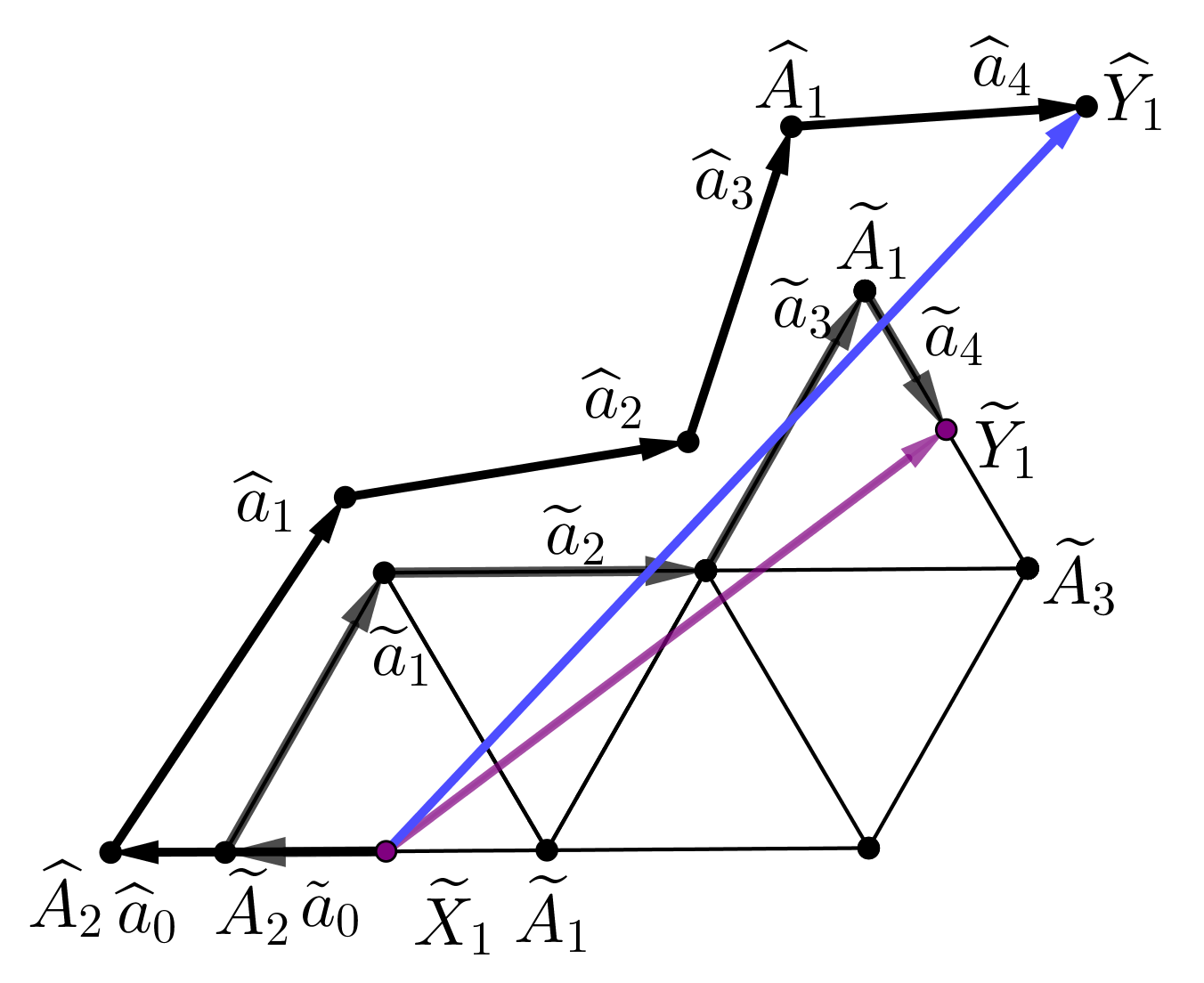}
\caption{ }
\label{dev_compare}
\end{center}
\end{figure}

Let the development  $\widetilde{T}_{pq}$ be placed such that 
the point  $ \widetilde{X}_1$  coincides with $\widehat{X}_1 $  of  $\widehat{ T}_{pq}$, 
and the vector  $\widehat{ X}_1\widehat{ A}_2$ has the same direction
 with $\widetilde{X}_1\widetilde{ A}_2$.
Similarly to the above we have
\begin{equation}\label{tilde_X_1_tilde_Y_1} 
\widetilde{ X}_1\widetilde{Y}_1 = \widetilde{a_0}+\widetilde{a_1}
+\dots+\widetilde{a_s}+\widetilde{a}_{s+1},
\end{equation}
where  $\widetilde{a_i}$ are  the sequential vectors of the $\widetilde T_{pq}$  boundary,  
   $s=\left[ \frac{p+q}{2} \right]+1$ and 
$\widetilde{a_0}~=~\widetilde{ X}_1 \widetilde{A}_2$, 
$\widetilde{a}_{s+1}~=~\widetilde{ A}_1 \widetilde{Y}_1$  (see Figure~\ref{dev_compare}).

Suppose the minimal distance from the vertices of $\widetilde{  T}_{pq}$  to the segment
$\widetilde{ X}_1 \widetilde{ Y}_1$ is reached
 at the vertex $ \widetilde{ A}_k$ and equals $\widetilde{ h}$
from the formula  (\ref{evcl_dist_vertex}).
Let us  estimate the distance $\widehat h$  between the segment 
  $\widehat { X}_1 \widehat {Y }_1$ 
and the corresponding vertex  $\widehat  {A }_k$ on  $\widehat{ T}_{pq}$.
 A geodesic on a regular tetrahedron in $\mathbb{E}^3$
  intersects at most three edges starting from the same
tetrahedron's vertex.
It follows, that the interior angles of the polygon $ \widetilde{ T}_{pq}$ are not greater than $4 \pi/3$. 
Hence the angles of the corresponding vertices 
on $\widehat{T}_{pq}$  are not   greater than $4 \widehat{\alpha}_i$. 
Applying  (\ref{bar_alpha-pi_3}) for  $1\le i \le s$
we get  that the angle between $ \widehat{a_i}$ and $ \widetilde{a_i} $ satisfies the inequality 
 \begin{equation}\label{angle_b_ia_i_1}
\angle(\widehat{a_i}, \widetilde {a_i} ) < \sum_{j=0}^i 4 \left( \pi  \tan^2 \frac{j}{R} + \varepsilon \right).  
 \end{equation} 
Since $R=1/a$, then using  (\ref{a_above_epsilon})  we obtain  
 \begin{equation}\label{tg_r_R}
  \tan \frac{j}{R} <  \tan \left( j \pi  \sqrt{2  \cos \frac{\pi}{12} } \; \sqrt{  \varepsilon  }   \right).   
  \end{equation}
The inequality (\ref{a_p+q}) holds if  the following condition fulfills
   \begin{equation}\label{vspomagat_3}
  \tan \left( j \pi  \sqrt{2  \cos \frac{\pi}{12} } \; \sqrt{  \varepsilon  }   \right) <   \tan \frac{\pi j}{2 (p+q) }.   
  \end{equation}
If $\tan x < \tan x_0$, then $\tan x < \frac{\tan x_0}{x_0} x$. 
From  (\ref{vspomagat_3}) it follows 
 \begin{equation}\label{vspomagat_4}
\tan \left( j \pi  \sqrt{2  \cos \frac{\pi}{12} } \; \sqrt{  \varepsilon  }   \right) <  
 2 (p+q)  \tan   \frac{\pi j}{2 (p+q) }  \sqrt{2 \cos \frac{\pi}{12}} \;  \sqrt{\varepsilon}.  
\end{equation}
Therefore from (\ref{tg_r_R}) and  (\ref{vspomagat_4}) we get  
 \begin{equation}\label{vspomagat_5}
\tan \frac{j}{R} < 2 (p+q)  \tan   \frac{\pi j}{2 (p+q) }  \sqrt{2 \cos \frac{\pi}{12}} \;  \sqrt{\varepsilon}.  
\end{equation}
Using (\ref{angle_b_ia_i_1}) and (\ref{vspomagat_5})
 we obtain the final estimation for the angle between 
the vectors $\widehat{a_i}$ and $ \widetilde{a_i}$:
\begin{equation}\label{angle_b_ia_i}
\angle(\widehat{a_i}, \widetilde {a_i} ) < \sum_{j=0}^i  4
  \left(   8 \pi     (p+q)^2 \cos \frac{\pi}{12}  \tan^2 \frac{\pi j}{2 (p+q) }  + 1 \right)  \varepsilon.       
  \end{equation}

Now estimate the  length of the vector $ \widehat{a_i} - \widetilde{a_i} $. 
The following inequality holds
 \begin{equation}\label{a_i-a_i_1}
| \widehat{a_i} - \widetilde{a_i} | \le 
\left|  \frac{\widehat{a_i}}{|\widehat{a_i}|} - \widetilde{a_i}\right|+ 
\left| \widehat{a_i} - \frac{\widehat{a_i}}{|\widehat{a_i}|} \right|.
 \end{equation} 
Since $\widetilde{a_i}$ is a unite vector, then
\begin{equation}
\left|  \frac{\widehat{a_i}}{|\widehat{a_i}|} - \widetilde{a_i}\right| 
\le \angle(\widehat{a_i}, \widetilde {a_i} ) 
 \;\;\; \textnormal{and}  \;\;\;
\left|\widehat{a_i} - \frac{\widehat{a_i}}{|\widehat{a_i}|}\right| \le \widehat{ l}_i -1.   
\end{equation} 
From the inequality (\ref{bar_l-1}) we get
  \begin{equation}\label{vspom_a_i_1}
\left|\widehat{a_i} - \frac{\widehat{a_i}}{|\widehat{a_i}|}\right|  < 
\frac {  \cos \frac{\pi}{12}  \left( 4+  \pi^2  (2i+1)^2  \right) }
{  \left( 1-\frac{2}{\pi}a(i+1) \right)^2 } \cdot \varepsilon.  
  \end{equation}
 Estimate the denominator in  (\ref{vspom_a_i_1}) using  (\ref{a_p+q}). Then
 \begin{equation}\label{vspom_a_i}
\left|\widehat{a_i} - \frac{\widehat{a_i}}{|\widehat{a_i}|}\right| <
\frac { \cos \frac{\pi}{12} (p+q)^2  \left( 4+  \pi^2  (2i+1)^2  \right) }{ \left( p+q-i-1 \right)^2 } \cdot \varepsilon.  
 \end{equation}
From (\ref{a_i-a_i_1}), (\ref{angle_b_ia_i}) and (\ref{vspom_a_i}) we obtain
  \begin{equation}\label{a_i-a_i}
| \widehat{a_i} - \widetilde{a_i} | \le \left( c_l (i) + \sum_{j=0}^i  c_\alpha (j) \right) \varepsilon,
 \end{equation}
where
   \begin{equation}\label{c_l}
 c_l (i) =\frac { \cos \frac{\pi}{12} (p+q)^2  \left( 4+  \pi^2  (2i+1)^2  \right) }{ \left( p+q-i-1 \right)^2 },
 \end{equation}
    \begin{equation}\label{c_alpha}
 c_\alpha (j)  =4   \left(   8 \pi     (p+q)^2 \cos \frac{\pi}{12}  \tan^2 \frac{\pi j}{2 (p+q) }  + 1 \right).
 \end{equation}
We estimate the length of $\widehat{ Y}_1 \widetilde{ Y}_1 $  using  (\ref{a_i-a_i}) 
 \begin{equation}\label{Y_1Y_1} 
 | \widehat{ Y}_1\widetilde{Y}_1| < 
  \sum_{i=0}^{s+1}  |\widehat{a_i} - \tilde {a_i} | <
 \sum_{i=0}^{s+1}   \left( c_l (i) + \sum_{j=0}^i  c_\alpha (j) \right) \varepsilon. 
 \end{equation}
From  (\ref{angle_b_ia_i}) it follows that the angle $\angle  \widehat{Y}_1\widehat{X }_1\widetilde{Y}_1$  satisfies
 \begin{equation}\label{angle_Y_1X_1Y_1} 
  \angle \widehat{Y}_1\widehat{ X}_1 \widetilde Y_1 < 
 \sum_{i=0}^{s+1}   c_\alpha (i)  \varepsilon.
 \end{equation}
The distance between the vertices $\widehat{A}_k$ and   $\widetilde{ A}_k$  equals
  \begin{equation}\label{A_kA_k} 
 |\widehat{ A}_k\widetilde{A}_k|<
 \sum_{i=0}^k  \left( c_l (i) + \sum_{j=0}^i  c_\alpha (j) \right) \varepsilon. 
 \end{equation}

We drop a perpendicular $\widehat{ A}_k\widehat{ H}$ from the  vertex $\widehat{ A}_k$ 
into the segment  $\widehat{X}_1\widehat{ Y}_1$.
The length of $\widehat{A}_k\widehat{H}$  equals $\widehat{ h }$.
Then we drop the perpendicular $\widetilde{ A}_k\widetilde{ H}$ 
into the segment $\widetilde{X}_1\widetilde{ Y}_1$ 
and the length of $\widetilde{ A}_k\widetilde{ H}$ equals $\widetilde{ h}$
  (see Figure~\ref{estim_dist}). 


 \begin{figure}[h]
\begin{center}
\includegraphics[width=85mm]{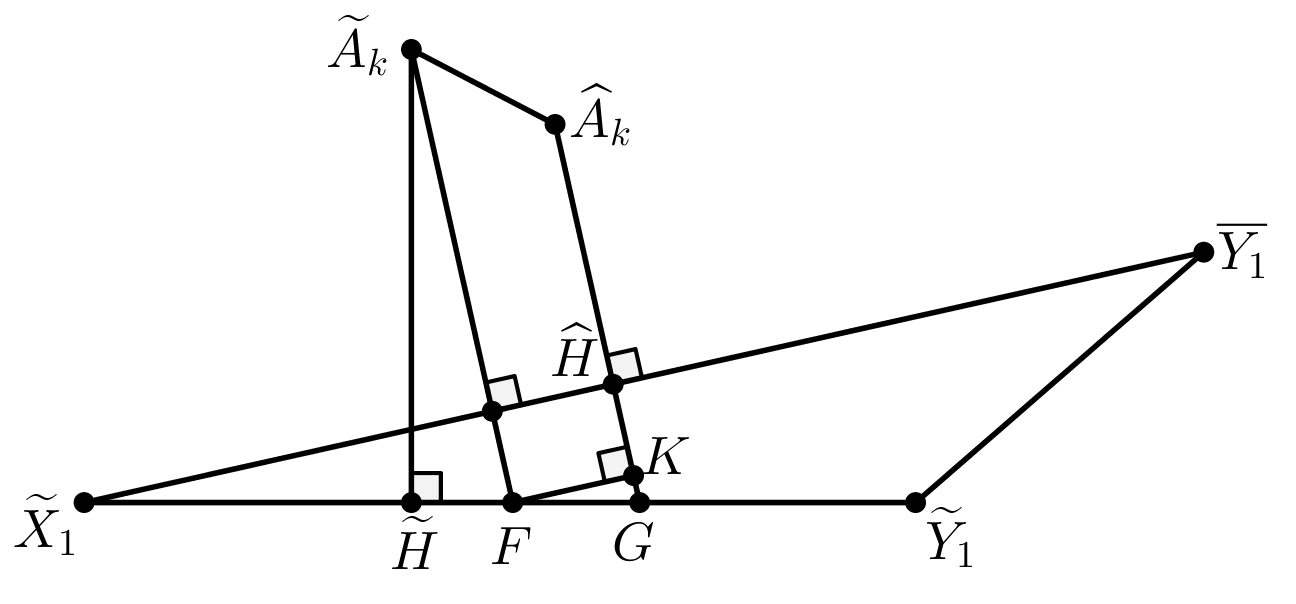}
\caption{ }
\label{estim_dist}
\end{center}
\end{figure}

Let the point $F$ on $\widetilde{ X}_1\widetilde{Y}_1$ be such that
the segment $\widetilde{ A}_k F$ is perpendicular to $\widehat{ X}_1\widehat{ Y}_1$.
Then the length of $\widetilde{ A}_k F$ is at least ~$\widetilde{ h }$.
Let   $G$ be the point of intersection  of
 $\widetilde{X}_1\widetilde{ Y}_1$ and  the extension of  $\widehat{ A}_k\widehat{ H}$.
 Let $FK$ is perpendicular to $\widehat{ H} G$  (see Figure~\ref{estim_dist}).
Then the length of $FK$ is not greater than the length of $ \widehat{A}_k\widetilde{A}_k $, 
 and $\angle KFG = \angle\widehat{ Y}_1\widehat{X}_1 \widetilde{ Y}_1$.
From the triangle $GFK$ we obtain
   \begin{equation}\label{FG_1} 
 | FG |= \frac{|FK|}{\cos \angle \widehat{ Y }_1\widehat{ X}_1 \widetilde{ Y}_1  }. 
 \end{equation}
Applying the inequality  $\cos x > 1-\frac{2}{\pi}x$ for $x < \frac{\pi}{2}$,  to (\ref{FG_1}),  we obtain
    \begin{equation}\label{FG_2} 
 | FG | < \frac{|\widehat{ A}_k\widetilde{ A}_k |}{ 1-\frac{2}{\pi}  
\angle \widehat{ Y}_1\widehat{ X}_1 \widetilde{ Y }_1 }. 
 \end{equation}
The inequalities  (\ref{angle_Y_1X_1Y_1}),  (\ref{A_kA_k})  and (\ref{FG_2})  imply 
   \begin{equation}\label{FG_3} 
   | FG | <  \frac{  \sum_{i=0}^k  \left( c_l (i) + \sum_{j=0}^i  c_\alpha (j) \right) \varepsilon}
  {1- \sum_{i=0}^s 
   \left(   64 \pi     (p+q)^2 \cos \frac{\pi}{12}  \tan^2 \frac{\pi i}{2 (p+q) }  + \frac{8}{\pi} \right) \varepsilon}. 
  \end{equation}
Applying  (\ref{epsilon_estim_1}) to the denominator in (\ref{FG_3}) we obtain 
       \begin{equation}\label{FG}
   | FG |<  \frac{   \sum_{i=0}^k  \left( c_l (i) + \sum_{j=0}^i  c_\alpha (j) \right) \varepsilon }
  {1-  \frac{(p+q+2)}{2\pi  \cos \frac{\pi}{12}  (p+q)^2 } -8 \sum_{i=0}^{s+1}  \tan^2 \left( \frac{\pi i}{2 (p+q) } \right)  }. 
  \end{equation}
 Therefore we have
   \begin{equation}\label{vspom_h} 
 \widetilde{ h} \le \widetilde{A}_kF \le  \widehat{ h} + |\widehat{ H} G |+
 |\widehat{ A}_k \widetilde {A}_k| + | FG |;  
 \end{equation}
Note, that $|\widehat{ H} G |<   | \widehat{Y}_1\widetilde{ Y}_1 |  $. 
 Lemma~\ref{evcl_dist_vertex_lem} implies that 
\begin{equation}\label{ }
\widetilde{ h} >   \frac{ \sqrt{3}  }{  4 \sqrt{p^2+q^2+pq} }.  \notag
\end{equation}
From (\ref{vspom_h}) it follows, that
   \begin{equation}\label{bar_h_1} 
      \widehat{ h }    > \frac{ \sqrt{3}  }{  4 \sqrt{p^2+q^2+pq} } -
 | \widehat{Y}_1\widetilde{ Y}_1 |   - 
|\widehat{ A}_k \widetilde {A}_k| - | FG |.  
 \end{equation}
Applying the estimations (\ref{Y_1Y_1}), (\ref{A_kA_k}), (\ref{FG}) 
and the identity $s=\left[ \frac{p+q}{2} \right]+1 $, we obtain
   \begin{equation}\label{bar_h_1} 
      \widehat{ h }    > \frac{\sqrt{3}}{4 \sqrt{p^2+q^2+pq}  } - c_0
 \sum_{i=0}^{\left[ \frac{p+q}{2} \right]+2}   \left( c_l (i) + \sum_{j=0}^i  c_\alpha (j) \right) \varepsilon,  
 \end{equation}
where  $c_l (i)$ is from  (\ref{c_l}), and $c_\alpha (j)$ is from (\ref{c_alpha}) and
          \begin{equation}\label{c_0 }
 c_0=   \frac{ 
 3 - \frac{(p+q+2)}{\pi  \cos \frac{\pi}{12}  (p+q)^2 } - 
   16 \sum_{i=0}^{\left[ \frac{p+q}{2} \right]+2}  \tan^2 \left( \frac{\pi i}{2 (p+q) } \right) }
  {1-  \frac{(p+q+2)}{2\pi  \cos \frac{\pi}{12}  (p+q)^2 } -
  8 \sum_{i=0}^{\left[ \frac{p+q}{2} \right]+2}  \tan^2 \left( \frac{\pi i}{2 (p+q) } \right)  }, \notag
 \end{equation} 
The inequality (\ref{bar_h_1}) implies that if $\varepsilon$ satisfies the condition
  \begin{equation}\label{varepsilon} 
     \varepsilon  <   \frac{\sqrt{3}}{4  c_0\sqrt{p^2+q^2+pq}\;
 \sum_{i=0}^{\left[ \frac{p+q}{2} \right]+2}   \left( c_l (i) + \sum_{j=0}^i  c_\alpha (j) \right)} ,   
 \end{equation} 
  then the distance from the vertices of the polygon
    $\widehat{ T}_{pq}$ to  $\widehat{ X}_1\widehat{ Y}_1$ is nonzero.

     Since we use the estimation  (\ref{epsilon_estim_1}), we get, that if 
    \begin{equation}\label{varepsilon_fin} 
     \varepsilon  <  \min \left\{
      \frac{\sqrt{3}}{4 c_0\sqrt{p^2+q^2+pq}\;
 \sum_{i=0}^{\left[ \frac{p+q}{2} \right]+2}   \left( c_l (i) + \sum_{j=0}^i  c_\alpha (j) \right)};
    \frac{1}{8  \cos \frac{\pi}{12}(p+q)^2 }
  \right\},   
 \end{equation} 
 then the segment  $\widehat{ X}_1\widehat{ Y}_1$ lies inside the polygon  $\widehat{ T}_{pq}$.
This implies  that the arc $X_1Y_1$  on a sphere lies inside the polygon $T_{pq}$.
 The arc $X_1Y_1$ corresponds to a simple closed geodesic  $\gamma$ of type $(p,q)$  
 on a regular tetrahedron with the planar angle $\alpha  =\pi/3 + \varepsilon$ in spherical space.
 From Corollary \ref{uniqueness_spher} we get, that this geodesic is unique, 
up to the rigid motion of the tetrahedron.
 
 Note, that the geodesic  $\gamma$ is invariant under 
  the rotation of the tetrahedron of  the angle $\pi$ over the line passing through 
 the midpoints of the opposite edges of the tetrahedron.
 The rotation of the tetrahedron through the angle  $2\pi/3 $ or $4\pi/3 $  over the altitude
dropped from the  vertex  to the center of its opposite face changes  $\gamma$ into
 another simple closed geodesic of  type $(p, q)$.
 
Rotating over the lines connecting other  vertices of the tetrahedron
 with the center of the opposite faces
doesn't give us any new geodesic. 
So if  $\varepsilon$ satisfies the condition (\ref{varepsilon_fin}), then
on a regular tetrahedron with the planar angle   $\alpha~=~\pi/3~+~\varepsilon$  in a spherical space
there exist three different simple closed geodesics of type  $(p, q)$,
disregarding isometries of the tetrahedron.
    \end{proof}

 \subsection{Necessary and sufficient condition for the existence of a simple closed geodesic.} 
 Let  $T(\alpha)$ be a regular tetrahedron with planar angles $\alpha$
  in spherical space~$\mathbb{S}^3$ of curvature~$1$.
Consider a development $R_{p,q}(\alpha)$ of  $T(\alpha)$ in $\mathbb{S}^3$ along 
a simple closed geodesic $\gamma_{p,q}$ of type $(p,q)$,
for  $\alpha\in(\pi/3, \pi/3+ \varepsilon)$,
 where $ \varepsilon$ is from Theorem~\ref{sufficient_cond}.
It follows from  Lemma~\ref{sym_develop} that the development $R_{p,q}(\alpha)$ has
 four points of symmetry $X_1(\alpha)$, $X_2(\alpha)$, $Y_1(\alpha)$, $Y_2(\alpha)$ and $X'_1(\alpha)$
 that correspond  to the midpoints of
 two pairs of opposite edges of the tetrahedron.
The geodesic  $\gamma_{p,q}$  passes through these midpoints.

\begin{figure}[h]
\begin{center}
\includegraphics[width=55mm]{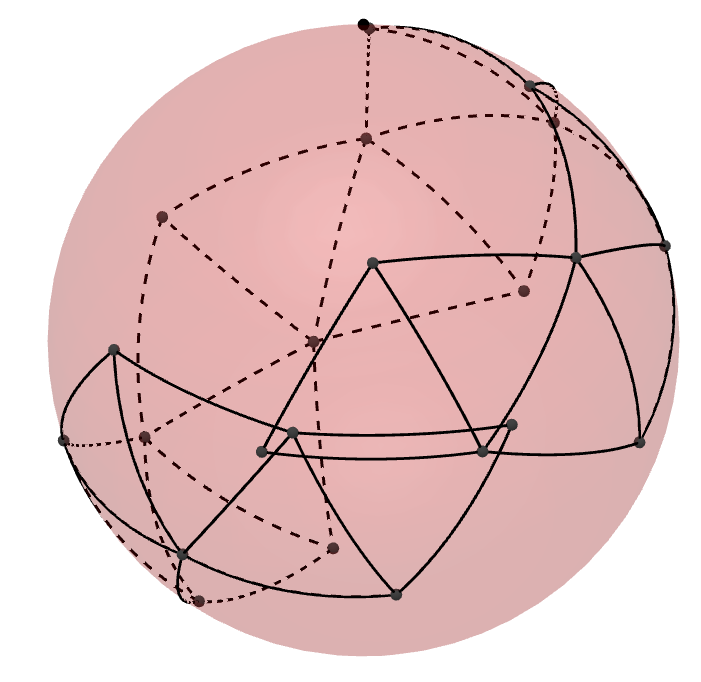}
\caption{ }
\label{general_devel}
\end{center}
\end{figure}

Now for fixed $(p,q)$ consider a one-parameter family of closed polygons 
 $R_{p,q}(\alpha)$, where $\alpha \in (\pi/3, 2\pi/3)$. 
Then $R_{p,q}(\alpha)$  may have overlaps on the sphere. 
However $R_{p,q}(\alpha)$ is considered as an abstract polygon homeomorphic to a disc, with
intrinsic metric, since each interior point of this polygon has a neighbourhood isometric
to the interior of a disc on the unit sphere $\mathbb{S}^2$.
This polygon is locally isometrically immersed in the sphere $\mathbb{S}^2$
   (see Figure~\ref{general_devel}).  
The development  $R_{p,q}(\alpha)$ also has a symmetry property for any $\alpha \in (\pi/3, 2\pi/3)$
with corresponding points 
$X_1(\alpha)$, $X_2(\alpha)$, $Y_1(\alpha)$, $Y_2(\alpha)$ and $X'_1(\alpha)$ on them. 

Then consider rectifiable curves $\sigma_{p,q}(\alpha)$ on $R_{p,q}(\alpha)$ that connect the
points  $X_1(\alpha)$, $X'_1(\alpha)$ and 
pass through  $X_2(\alpha)$, $Y_1(\alpha)$,  $Y_2(\alpha)$.
If   $X_1(\alpha)X'_1(\alpha)$ lies inside the 
 development $R_{p,q}(\alpha)$, then $\sigma_{p,q}(\alpha)$ corresponds to the
 simple closed geodesic on regular tetrahedron $T(\alpha)$. 
  From Theorem~\ref{sufficient_cond} follows that this is true if $\alpha $ is   close to $\pi/3$. 
 Then from  Lemma~\ref{length} we get that length of $\sigma_{p,q}(\alpha)$  is less then $2\pi$. 
In~\cite{Bor2022} Borisenko proved that this condition is also sufficient for existence
 a simple closed geodesic on a regular tetrahedron in $\mathbb{S}^3$.

 The infimum $L_{p,q}(\alpha)$ of the lengths of the curves $\sigma_{p,q}(\alpha)$ 
is referred to as \textit{the length of the abstract shortest curve} in the development.

%

\begin{theorem}\label{ness_suff_cond}\textnormal{\cite{Bor2022}}
On a regular tetrahedron in spherical space of curvature one there exist 
a simple closed geodesic of type $(p,q)$
if and only if the length of the abstract shortest curve on the development is less than $2\pi$.
\end{theorem}

 \begin{proof}
 1. \textit{Necessity.} 
If there exist a simple closed geodesic of type $(p,q)$ on a tetrahedron $T(\alpha)$, then 
by unfolding along this geodesic we obtain $R_{p,q}(\alpha)$. 
The geodesic unfolds into an arc of great circle, which lies inside $R_{p,q}(\alpha)$, 
connects the points $X_1(\alpha)$ and $X'_1(\alpha)$
 and passes through the points of symmetry of $R_{p,q}(\alpha)$.  
  Lemma~\ref{length} implies that 
   $L_{p,q}(\alpha)$ equals the length of this geodesic 
and  $L_{p,q}(\alpha)$  is less then $2\pi$ 
  (see Figure~\ref{geod(2,3)_spher-2}). 
 \begin{figure}[h]
\begin{center}
\includegraphics[width=55mm]{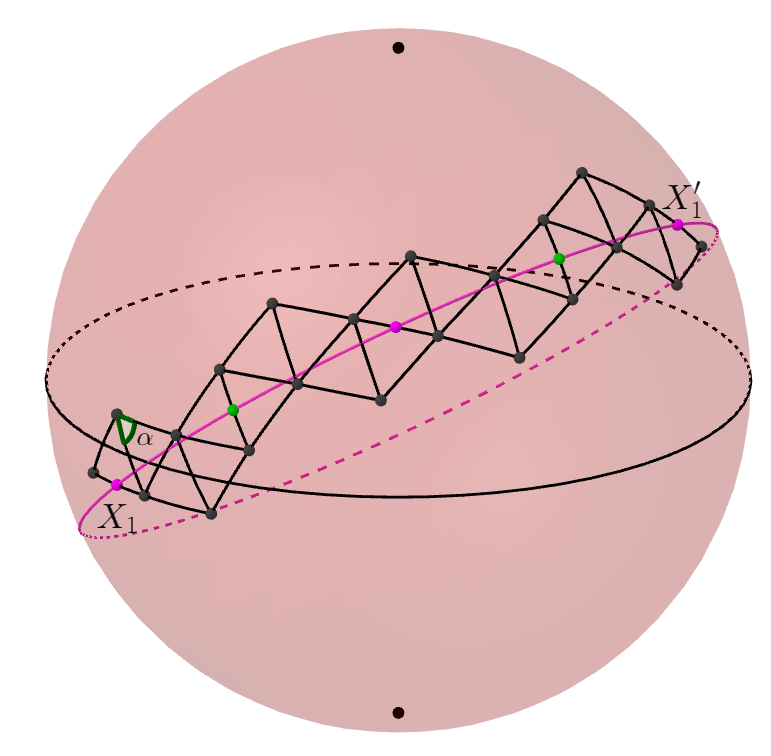}
\caption{ }
\label{geod(2,3)_spher-2}
\end{center}
\end{figure}

2. \textit{Sufficiency.} 
Let us proof the monotonicity of $L_{p,q}(\alpha)$.
Let the infimum   $L_{p,q}(\alpha)$ is attained on a curve $\sigma_{p,q}(\alpha)$  on $R_{p,q}(\alpha)$.
Consider the geodesic mapping of the sphere $\mathbb{S}^3$ onto the 
Euclidean tangent space $T_{O}\mathbb{S}^3$, where $O$ is a center of the
inscribed sphere in the tetrahedron  $T(\alpha)$. 
Then $T(\alpha)$ is mapped onto the regular tetrahedron  $\widehat T(\alpha)$ in $\mathbb{E}^3$ 
and the curve $\sigma_{p,q}(\alpha)$ is mapped onto $\widehat\sigma_{p,q}(\alpha)$.

Let $\widehat T(\alpha(\lambda))=\lambda\widehat T(\alpha)$ be the tetrahedron
homothetic to  $T(\alpha)$ with center $O$ and ratio $\lambda < 1$, so that $\alpha(\lambda)< \alpha$. 
This homothety takes $\widehat\sigma_{p,q}(\alpha)$ to
 a curve $\widehat\sigma_{p,q}(\alpha(\lambda))$. 
 
 Consider the inverse geodesic mapping of  $T_{O}\mathbb{S}^3$ onto  $\mathbb{S}^3$.
 It takes $\widehat T(\alpha(\lambda))$ to a regular tetrahedron
  $T(\alpha(\lambda))$ where $\alpha(\lambda)< \alpha$.
  The curve $\widehat\sigma_{p,q}(\alpha(\lambda))$ is mapped to $\sigma_{p,q}(\alpha(\lambda))$ 
  that belongs to our class of curves. 
  Let us show that the length of the curve $\sigma_{p,q}(\alpha(\lambda))$ is less than 
  $L_{p,q}(\alpha)$ for $\lambda < 1$.
 
 The curve $\widehat\sigma_{p,q}(\alpha)$ consists of the finite number of segments with 
 endpoints on edges of the regular tetrahedron. 
Consider one of these segments  $\widehat z(\alpha)$ on
 the face $A_1A_2A_3$ of $\widehat T(\alpha)$. 
The family of segments $\lambda\widehat z (\alpha)$ on $\lambda\widehat T(\alpha)$ 
is  homothetic to $\widehat z(\alpha)$ with respect to the center $O$. 
The great circle arc  $z(\lambda) = z(\alpha(\lambda))$ is the inverse geodesic
 images of $\lambda\widehat z(\alpha)$.
We show that the length of $z(\lambda)$ is monotonically increasing function of~$\lambda$.

 \begin{figure}[h]
\begin{center}
\includegraphics[width=55mm]{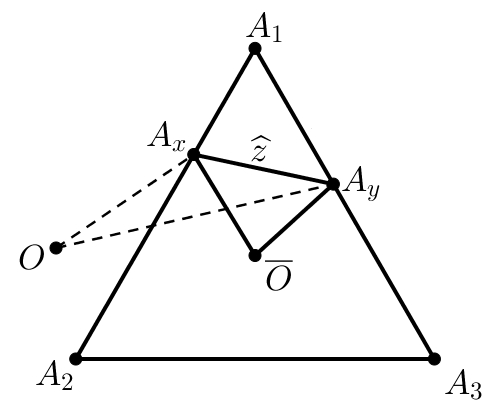}
\caption{ }
\label{homotet}
\end{center}
\end{figure}

Denote by $ A_x$ and $ A_y$ the endpoints of  $\widehat z(\alpha)$
on  $A_1A_2$ and $A_1A_3$ respectively. 
Then 
$$ | A_x  A_y|^2 =  | A_1 A_x|^2 + 
 |  A_1 A_y|^2 - | A_1 A_x|| A_1  A_y|.$$
  The radius of the inscribed sphere of the tetrahedron $\widehat T(\alpha)$
  with edge length $a$ equals $r=a/(2\sqrt{6})$.
  The distance from the center of $\widehat T(\alpha)$ to the points
    $  A_x$ and $  A_y$  can be found from the triangles 
    $\triangle A_1\bar O A_x$, where
    $\bar O$ is the center of the face $A_1A_2A_3$   (see Figure~\ref{homotet}):
    $$ |\bar O  A_x|^2 =| A_1 A_x|^2 +
     \frac{a^2}{3} - a | A_1 A_x|.$$
  From the triangle   $\triangle O\bar O A_x$ we get
   $$ |\bar O  A_x|^2 =\frac{3}{8} a^2 +|  A_1  A_x|^2  -  a | A_1  A_x|.$$
    From the triangle   $\triangle O\bar O A_y$ we have 
   $$ |\bar O A_y|^2 =\frac{3}{8} a^2 +| A_1 A_y|^2  -
   a | A_1 A_y|.$$
  From the triangles $\triangle OS A_x$ and $\triangle OS A_y$, where
  $S$ is the center of the sphere $\mathbb{S}^3$, we obtain
  $$ |S A_x|^2=1+ |\bar O   A_x|^2; \;\; 
  |S A_y|^2=1+ |\bar O  A_y|^2. $$
  From $\triangle  A_x S  A_y$ we obtain 
  $$\cos z = \frac{\left( 1+ |\bar O   A_x|^2 \right)
  +\left( 1+ |\bar O A_y|^2 \right) - | A_x A_y|^2 }
  {2\sqrt{1+ |\bar O  A_x|^2 } \sqrt{1+ |\bar O A_y|^2 } },$$
  where $z$ is the angle at the vertex $S$. 
  
  Similarly, for the homothetic tetrahedra $\lambda \widehat T(\alpha)$, we have
   $$\cos z(\lambda) = \frac{\left( 1+ \lambda^2|\bar O \widehat A_x|^2 \right)
  +\left( 1+ \lambda^2|\bar O   A_y|^2 \right) 
  - \lambda^2|  A_x A_y|^2 }
  {2\sqrt{1+\lambda^2|\bar O  A_x|^2 }
   \sqrt{1+ \lambda^2|\bar O  A_y|^2 } }.$$
  
 The derivative of $z(\lambda)$ at $\lambda = 1$ is positive. 
 This implies that the length of $\sigma_{p,q}(\alpha(\lambda))$ is less than 
 the length of $\sigma_{p,q}(\alpha )$ for $\lambda < 1$.
 Hence $L_{p,q}(\alpha(\lambda))<L_{p,q}(\alpha)$ for $\lambda < 1$ and $\alpha(\lambda)<\alpha$.
 
For $\pi/3<\alpha<\pi/3+\varepsilon$, 
 where  $\varepsilon$ is from Theorem~\ref{sufficient_cond},
 there is a simple closed geodesic of type $(p,q)$ on a regular tetrahedron in $\mathbb{S}^3$.
 This geodesic unfolds into a  curve $\sigma_{p,q}(\alpha)$ 
 of length  $L_{p,q}(\alpha)<2\pi$  inside the development $R_{p,q}(\alpha)$.
 
Now increase the angle $\alpha$ starting from  $\alpha_1$.
 As long as $\sigma_{p,q}(\alpha)$  lies inside the development $R_{p,q}(\alpha)$,
it corresponds to a simple closed geodesic on a regular tetrahedron $T(\alpha)$. 
Let $\beta$  be the first value of $\alpha$  for which  $\sigma_{p,q}(\alpha)$
attains the boundary of $R_{p,q}(\alpha)$.
This value exists by Theorem~\ref{necessary_cond} which implies 
 that there is $\alpha_2\in (\pi/3, \pi/2)$
 such that there is no simple closed geodesics on $T(\alpha)$ for $\alpha~>~\alpha_2$.
 
 The point of intersection of $\sigma_{p,q}(\beta)$ with the boundary of the development 
  $R_{p,q}(\beta)$  is a vertex of the tetrahedron. 
 Since $R_{p,q}(\beta)$ consists of congruent polygons, then the segment $\sigma_{p,q}(\beta)$
 `touches' the boundary of $R_{p,q}(\beta)$ at four vertices.
The property of symmetry of  $R_{p,q}(\beta)$ implies that
these `touchings' alternate so there are two  of them from each side of $\sigma_{p,q}(\beta)$ 
   (see Figure~\ref{dev_spher_plane_intersec}). 
 \begin{figure}[h]
\begin{center}
\includegraphics[width=65mm]{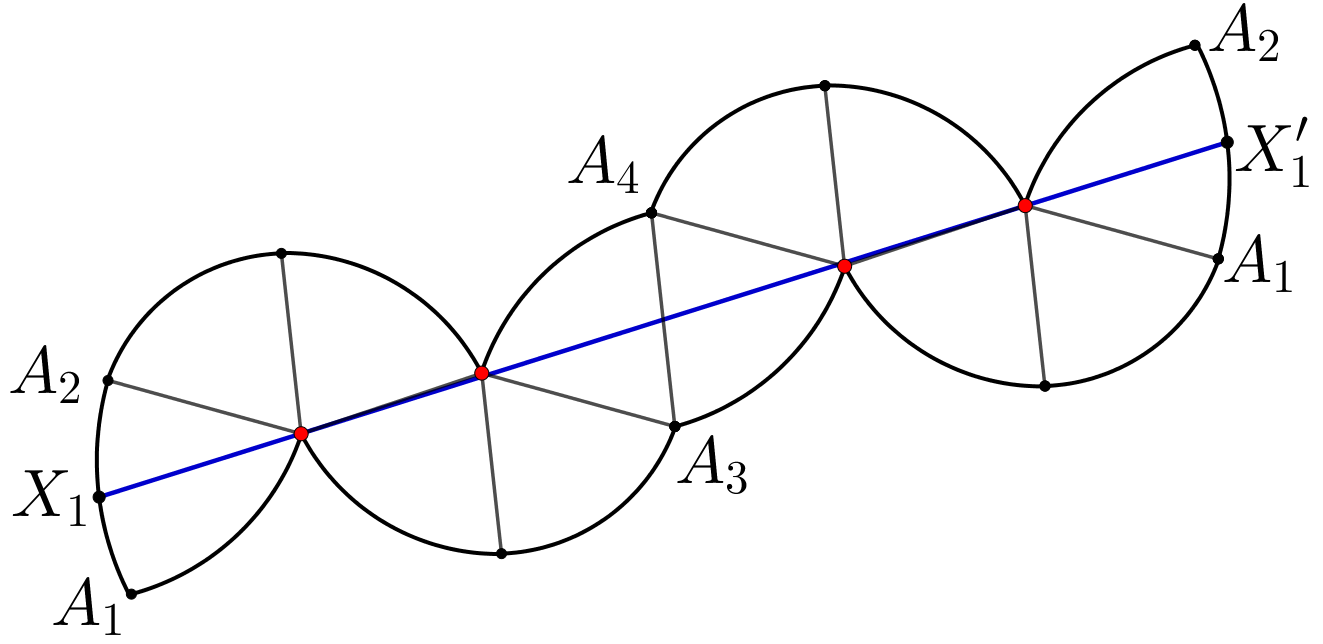}
\caption{ }
\label{dev_spher_plane_intersec}
\end{center}
\end{figure}

 The segment  $\sigma_{p,q}(\alpha)$ cannot `touch' the boundary of the development 
 $R_{p,q}(\beta)$ at five points.
 Otherwise the curve $\sigma_{p,q}(\alpha)$  passes twice through some vertex of $T(\beta)$.
  For any line segment, the full angle on one side is $\pi$. 
The full angle at any vertex is less than $2\pi$.
Thus the segments $l_1$ and $l_2$ of the curve $\sigma_{p,q}(\alpha)$
 intersect at a nonzero angle at that vertex. 
The geodesics $\sigma_{p,q}(\alpha)$ with
$\alpha<\beta$ and $\alpha$ close to $\beta$   also intersect themselves,
 which contradicts the fact that these geodesics are simple.

 The case when two points of intersection 
(for example, the vertices $A_2$ and $A_3$)
merge is also impossible.
 These two vertices are not connected by an edge, because if
we take $\alpha<\beta$ and let $\lim\alpha=\beta$, 
then we see that the length of the edge connecting
these two vertices of intersection tends to $0$.
As $\alpha~\rightarrow~\beta$, the full angles
at $A_2$ and $A_3$ tend to angles $\ge \pi$,
 for otherwise the geodesics $\sigma_{p,q}(\alpha)$ would
cross the boundary of the development for some  $\alpha<\beta$.
 Without loss of generality,
we can assume that $\beta\le \beta_0=2\arcsin \sqrt{7/18}<2\pi$,
 since there are only three simple
closed geodesics for $\beta \ge \pi/2$ (see  Lemma~\ref{basic_geod}). 
This bound follows from the case
$p=2, q=1$ of inequality~(\ref{estim_above}) from Theorem~\ref{necessary_cond}.
 For the full angles at the
vertices $A_2$   to tend to limits$\ge \pi$,
 it is necessary that at least three triangles
meet at $A_2$ and that for $\alpha$  close to $\beta$ 
 two edges meeting at  $A_2$ belong to
triangles in the development traversed by the line segment $\sigma_{p,q}(\alpha)$.
The same we observed for $A_3$. 
 Then four different edges of triangles would meet at the merged vertex.
Thus, four edges come out of a vertex of the tetrahedron, which is a contradiction.
 
As a result, for $\alpha = \beta$ the segment  $\sigma_{p,q}(\alpha)$ `touches'
the boundary of  $R_{p,q}(\beta)$ at four points,
which correspond to the vertices of the tetrahedron. 
The curve $\sigma_{p,q}(\alpha)$ divides the tetrahedron 
into two regions homeomorphic to a circle.
Each interior point has a neighborhood isometric to a disc on the sphere $\mathbb{S}^2$
of curvature $1$, and the boundary is a digon. 
The edges of this digon have the same
length, the full angles at both vertices are $3\beta - \pi$,
and the geodesic curvature of the digon is~$0$. 
Therefore, the perimeter of the digon is $2\pi$.
Hence the length of $\sigma_{p,q}(\alpha)$ is $2\pi$,
which implies that $L_{p,q}(\alpha)=2\pi$.

If a simple closed geodesic exists for a fixed $\alpha$,
 then $L_{p,q}(\alpha)$ is equal to the length
of this geodesic, and therefore it is   $< 2\pi$ for $\alpha<\beta$. 
If $\alpha> \beta$, then, due to
the monotonicity of $L_{p,q}(\alpha)$, 
the length of $L_{p,q}(\alpha)$ is greater than $2\pi$, and there are
no simple closed geodesics of type $(p, q)$ on the tetrahedron $T(\alpha)$.
\end{proof}

\begin{corollary}\label{sufficient_cond_edge} \textnormal{\cite{Bor2022} }
If the edge $a$ of a regular tetrahedron in spherical space satisfies the
inequality
\begin{equation}
a< 2\arcsin \frac{\pi}
{\sqrt{p^2+pq+q^2}+\sqrt{(p^2+pq+q^2)+2\pi^2} }
\end{equation}
then this tetrahedron has a simple closed geodesic of type $(p, q)$.
\end{corollary}
\begin{proof}
Let $O$ be the centre of the inscribed and circumscribed spheres of a regular
tetrahedron $T(\alpha)$ in  spherical space $\mathbb{S}^3$.

Consider a geodesic mapping of the
open hemisphere of $\mathbb{S}^3$  containing $T(\alpha)$
 onto the tangent space $T_O\mathbb{S}^3$. The
tetrahedron $T(\alpha)$ is mapped to a regular tetrahedron $\widehat T(\alpha)$ 
with center at $O$ in Euclidean space $T_O\mathbb{S}^3$.
The midpoints of the edges  are mapped to the midpoints. 
Let $\widehat a $  be the edge length of $\widehat T(\alpha)$.

Let $\widehat\gamma_{p,q}(\alpha)$ be a simple closed geodesic of type $(p, q)$
 that passes through the
midpoints of two pairs of opposite edges of $\widehat T(\alpha)$.
Then the length of $\widehat\gamma_{p,q}(\alpha)$ is equal 
\begin{equation}\label{L_pq}
\widehat L_{p,q}(\alpha)= 2 \widehat a \sqrt{p^2+pq+q^2}.   
\end{equation}
Take $\alpha$ such that $\widehat L_{p,q}(\alpha) < 2\pi$.
 The inverse image $\gamma_{p,q}(\alpha)$ of the geodesic $\widehat\gamma_{p,q}(\alpha)$
on $T(\alpha)$ has length less than $\widehat L_{p,q}(\alpha)$, 
and therefore less than $2\pi$. 
The curve $\gamma_{p,q}(\alpha)$
belongs to the class of admissible curves $\sigma_{p,q}(\alpha)$
 in the definition of $L_{p,q}(\alpha)$. 
 Therefore, $L_{p,q}(\alpha)<2\pi$, 
 and Theorem~\ref{ness_suff_cond} implies that there exists a simple closed geodesic
of type $(p, q)$ on  $T(\alpha)$. 
It remains to use the inequality
$$ 2 \widehat a \sqrt{p^2+pq+q^2}  < 2\pi $$
 to obtain a bound on $\alpha$, or, equivalently, on $a$. 
Formula~\ref{a} implies that 
$$ 2\sin (a/2)\cos (a/2) = 1. $$

We apply a geodesic mapping of the sphere $\mathbb{S}^3$ 
from its centre $S$ onto the tangent
space $T_O\mathbb{S}^3$.
 Consider the triangle $\triangle SOB$,
  where $B$ is the midpoint of $A_1A_2$.
   Let $\widehat  B$ be the image of $B$ under the geodesic mapping (Figure~\ref{triangle_SOBA1A2});
then $$|O\widehat  B| = \tan |OB|. $$
The edge $A_1A_2$ of the spherical triangle maps to the edge $\widehat A_1\widehat A_2$ 
of the regular tetrahedron in Euclidean space,
and  $\widehat A_1\widehat A_2$  is perpendicular to $O\widehat  B$. 
From the triangle $\triangle S\widehat A_1 \widehat B$ we obtain
\begin{equation}\label{a_2_vs} 
\frac{ \widehat a }{2} = |\widehat A_1 \widehat B| = 
|S\widehat  B|\tan\frac{a}{2} = \frac{\tan(a/2)}{\cos |OB|}. 
\end{equation}

 \begin{figure}[h]
\begin{center}
\includegraphics[width=95mm]{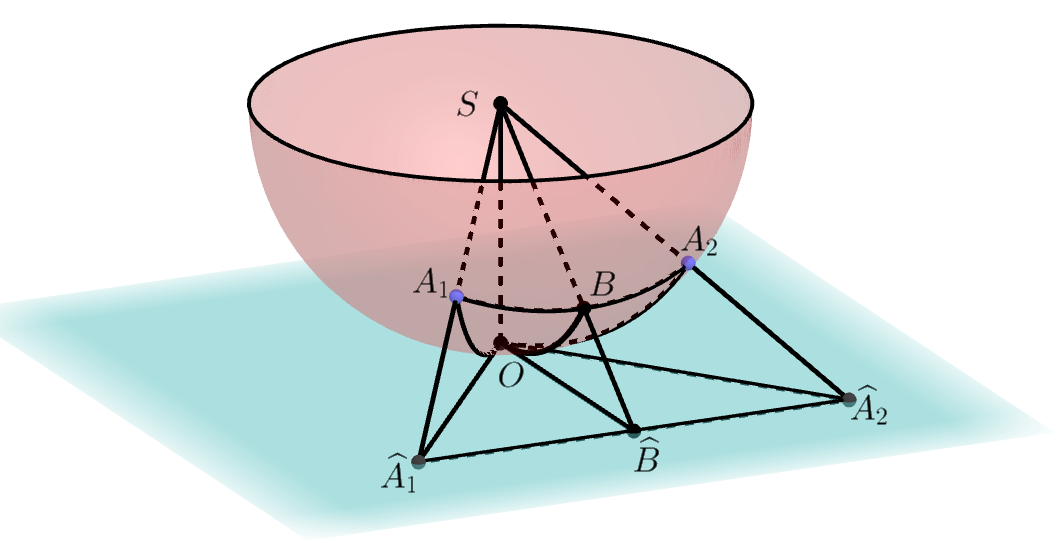}
\caption{ }
\label{triangle_SOBA1A2}
\end{center}
\end{figure}

From the triangle $\triangle PA_1A_2$ on a face of the tetrahedron in spherical space, 
where P is the centre of the inscribed and circumscribed circles of the
face, we obtain
$$ \cos a = \cos^2 R_{bas} - \frac{1}{2} \sin^2R_{bas}, $$
where $R_{bas}=|PA_1|=|PA_2|$.
 Hence
\begin{equation}\label{R_bas} 
  \cos R_{bas} =\sqrt{ \frac{1+2\cos a}{3} }.
    \end{equation}
From $\triangle A_4PA_1$ (Figure~\ref{triangle_A1A2A3}) we obtain
\begin{equation}\label{vspom}
\cos a = \cos(R+r)\cos  R_{bas},
\end{equation}
where $R$ is the radius of the circumscribed sphere of the tetrahedron $A_1A_2A_3A_4$,
$r$ is the radius of the inscribed ball, and $|A_4P| = R + r$.
Then~(\ref{vspom}) implies that
\begin{equation}\label{vspom_1}
 \cos R > \frac{ \cos a}{\cos  R_{bas}}.
\end{equation}
From $\triangle OA_1B$ we obtain
\begin{equation}\label{vspom_2}
 \cos R =\cos |OB| \cos (a/2).
\end{equation}
Expressions~(\ref{vspom_1}) and (\ref{vspom_2}) implies that 
\begin{equation}\label{vspom_3}
\frac{1}{\cos |OB|}= \frac{\cos (a/2)}{\cos R} < \frac{\cos(a/2) \cos R_{bas} }{ \cos a}.
\end{equation}
From~(\ref{a_2_vs}), (\ref{R_bas}) and (\ref{vspom_3}) we get
 \begin{equation}\label{vspom_4}
\widehat a /2 < \frac{\sin (a/2)}{\cos a} \sqrt{ \frac{1+2\cos a}{3} }
\le \frac{\sin (a/2)}{\cos a}.
\end{equation}
Therefore,  from~(\ref{L_pq}) and (\ref{vspom_4}) we obtain the following 
estimation for the length of a simple closed geodesic
$\widehat\gamma_{p,q}(\alpha)$ of type $(p,q)$ on $\widehat T(\alpha)$:
$$ \widehat L_{p,q}(\alpha) \le 4 \frac{\sin (a/2)}{\cos a}  \sqrt{p^2+pq+q^2}.  $$

 \begin{figure}[h]
\begin{center}
\includegraphics[width=75mm]{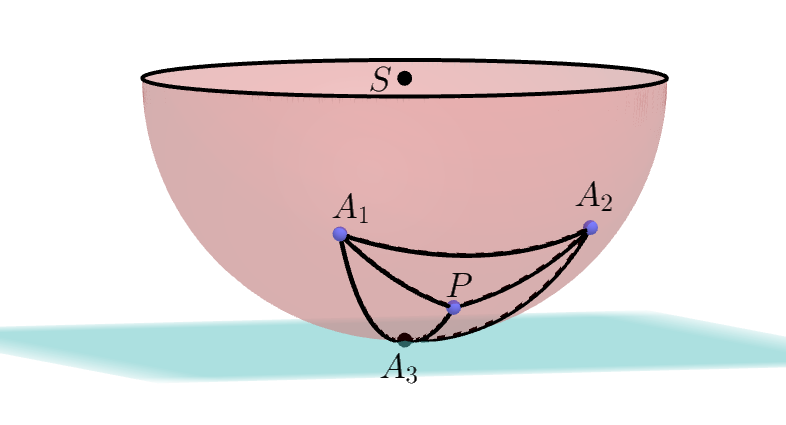}
\caption{ }
\label{triangle_A1A2A3}
\end{center}
\end{figure}

Remind, that from Theorem~\ref{ness_suff_cond} it follows that if
$ \widehat L_{p,q}(\alpha) < 2\pi$,  then there exist a simple closed geodesic 
of type $(p, q)$ on  $T(\alpha)$ in $\mathbb{S}^3$.
Resolving the quadratic inequality 
 $$  4 \frac{\sin (a/2)}{\cos a}  \sqrt{p^2+pq+q^2} < 2\pi  $$
with respect to $\sin(a/2)$, we obtain the required inequality. 
\end{proof}

\section{ Simple closed geodesics on regular tetrahedra in $\mathbb{H}^3$}
\subsection{Necessary conditions for a closed geodesic to be simple.}

We assume that the  Gaussian curvature of  {\it hyperbolic space} 
({\it Lobachevsky space}) $\mathbb{H}^3$ equals~$-1$.
A {\it  regular tetrahedron} in  $\mathbb{H}^3$
is a closed convex polyhedron  all of whose faces are regular 
 geodesic triangles  and all vertices are regular trihedral angles.
The planar angle $\alpha $ of the face  satisfies the inequality $0< \alpha<\pi/3$ and
the  length $a$ of   edges   is equal to
\begin{equation}\label{a_hyp}
a=\text{arcosh} \left(  \frac{\cos\alpha}{1-\cos\alpha}  \right).
\end{equation}

Consider  the  Cayley-Klein model of   hyperbolic space. 
In this model points are represented by the points in the interior of the unit ball.
Geodesics in this model are the chords of the ball.
Assume that the center  of the circumscribed sphere of a regular tetrahedron coincides 
with the center of the model.
Then the regular tetrahedron in hyperbolic space is represented by 
a regular tetrahedron in   Euclidean space.

\begin{lemma}\label{intersect_four_edges}\textnormal{\cite{BorSuh2020}}
If a geodesic on a regular tetrahedron in hyperbolic space intersects
three edges meeting at a common vertex consecutively, and intersects one of these
edges twice, then this geodesic has a point of self-intersection.
   \end{lemma}
\begin{proof}
Let $A_1A_2A_3A_4$ be a regular tetrahedron in  $\mathbb{H}^3$.
Suppose the geodesic $\gamma$   intersects $A_4A_1$, $A_4A_2$ and $A_4A_3$  consecutively
at the points $X_1$, $X_2$, $X_3$ respectively and 
then intersects the edge $A_4A_1$ again at the point $Y_1$.

Suppose that the length of $A_4X_1$ is  less than the  length of $A_4Y_1$.

Unfold the faces $A_1A_2A_4$, $A_4A_2A_3$ and $A_4A_3A_1$ to the hyperbolic plane.
 Consider the Cayley-Klein model of the hyperbolic plane and place the vertex $A_4$
 at the center of the model.
Then the part $X_1X_2X_3Y_1$ of the geodesic is the straight line segment on the development.
We obtain a triangle $X_1A_4Y_1$ on the development. 
 \begin{figure}[h]
\begin{center}
\includegraphics[width=50mm]{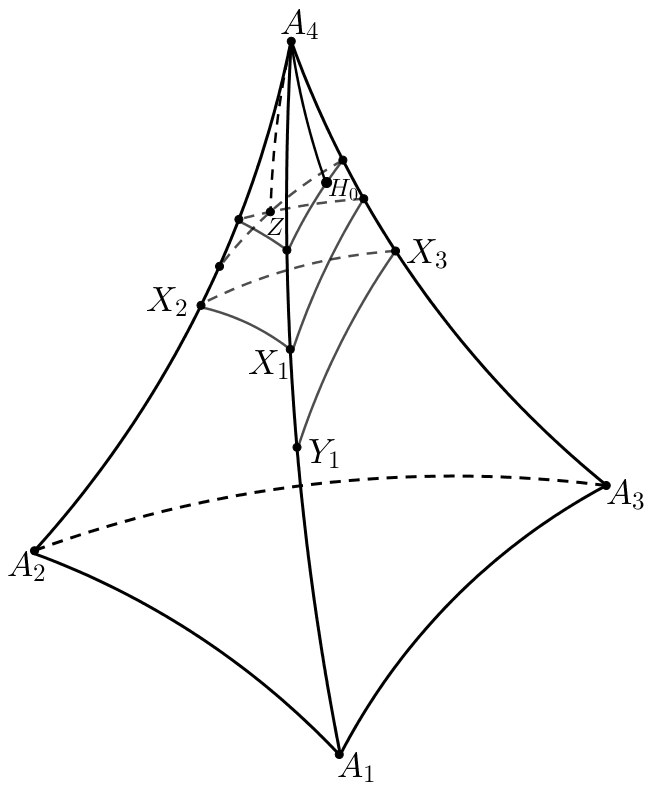}
\caption{ }
\label{X1X2X3Y1tetrV1H}
\end{center}
\end{figure}

Let $\rho(X)$ be the distance function between the vertex $A_4$ and a point $X$ on $\gamma$.
It is known that if   $\gamma$  be a  geodesic in a complete simply connected
 Riemannian manifold $M$ of nonpositive curvature, then the function $\rho(X)$ of a distance
from the fixed point $A$ on  $M$ to the points $X$ on $\gamma$ is a convex function.
The minimum of  $\rho(X)$ is achieved at the point $H_0$ 
such that $A_4H_0$ is orthogonal to $\gamma$ and 
$\angle H_0A_4Y_1 > 3\alpha/2$. 

 Let  $Z_1$  be the point on the segment  $H_0Y_1$  such that $\angle H_0A_4Z_1 = 3\alpha/2$.
On the opposite side of $H_0$  we choose the point  $Z_2$ such that 
$\angle H_0A_4Z_2 = 3\alpha/2$. 
The point $Z_2$ also lies on the face at the vertex $A_4$ of tetrahedron. 

Since $\angle H_0A_4Z_1 = \angle H_0A_4Z_2 = 3\alpha/2$, 
it follows that  the points $Z_1$ and $Z_2$ correspond to the same point  $Z$  on
the generatrix $A_4Z$ opposite to $A_4H_0$ on the tetrahedron.
This point is the self-intersection point of the geodesic $\gamma$ (Figure \ref{X1X2X3Y1tetrV1H}).
\end{proof}
 
\begin{lemma}\label{nesdist}\textnormal{\cite{BorSuh2020}}
Let $d$ be the minimum distance from the vertices of a regular tetrahedron
in hyperbolic space to a simple closed geodesic on the tetrahedron. Then
\begin{equation}\label{distvertex2}
d> \frac{1}{2} \ln \left( \frac{  \sqrt{2\pi^3} + \left( \pi- 3\alpha \right)^{\frac{3}{2}} }
{   \sqrt{2\pi^3} - \left( \pi- 3\alpha \right)^{\frac{3}{2}}    }   \right),
 \end{equation}
where $\alpha$ is the planar angle  of a face of the tetrahedron.
\end{lemma}
\begin{proof} 
Let $\gamma$ be a simple closed geodesic  
on a regular tetrahedron $A_1A_2A_4A_3$   in hyperbolic space $\mathbb{H}^{3}$.
Assume that minimum distance $d$ from the vertices of the tetrahedron to  $\gamma$ 
  is achieved at the vertex $A_4$ on the face  $A_2A_4A_3$.
Draw a generatrix $A_4H$   orthogonal to $\gamma$ at the point $H_0$.
 Denote by $\beta$ the angle  $\angle A_4HA_2$. 
Without loss of generality we assume that $0 \le \beta \le\alpha/2$. 
 \begin{figure}[h]
\begin{center}
\includegraphics[width=50mm]{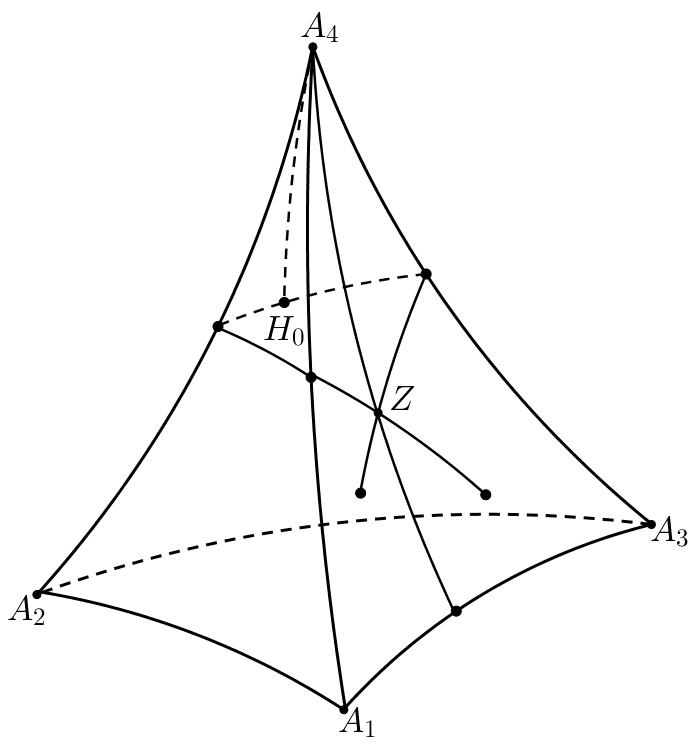}
\caption{ }
\label{nescondnew}
\end{center}
\end{figure}

We draw a generatrix  $A_4K$ such that
the planar angle between $A_4K$ and $A_4H$ equals  $3\alpha/2$.
Then   $A_4K$ lies  in the face $A_1A_4A_3$ and $\angle A_1A_4K = \alpha/2-~\beta$.
Note that if $\beta=\alpha/2$, then $A_4K$ coincides with $A_4A_1$.
If $\beta = 0$, then $A_4K$ coincides with the altitude in a face of the tetrahedron
 and has the smallest length~$h$
(Figure~\ref{nescondnew}). 

We cut the trihedral angle at  $A_4$ along the generatrix $A_4K$ and
 develop it to the hyperbolic plane in the Cayley-Klein model. 
 We put the vertex $A_4$  at the centre of the boundary circle. 
 The trihedral angle unfolds into a convex polygon $K_1A_4K_2A_3A_2A_1$.
The angle $K_1A_4K_2$ equals $3\alpha$. 
The segment $A_4H$ corresponds to the bisector of the angle $K_1A_4K_2$.
The geodesic $\gamma$ is a straight line  ortogonal to $A_4H$ at   $H_0$.

On the lines  $A_4K_1$ and $A_4K_2$ choose the points $P_1$ and $P_2$  respectively such that
 $|A_4P_1|=|A_4P_2| = h$.
The line segment  $P_1P_2$ is ortogonal to $A_4H$ at the point $H_p$, and
\begin{equation}\label{A4Hp}
\tanh |A_4H_p| = \cos(3\alpha/2)\tanh  h. \notag
\end{equation}

If $ d \le |A_4H_p|$, then $\gamma$ lies above the segment $P_1P_2$,
and therefore   $\gamma$ intersects the  lines $A_4K_1$ and $A_4K_2$ 
at the points $Z_1$ and $Z_2$ respectively.
When we  fold the development back to the tetrahedron, 
the  segments $A_4K_1$ and $A_4K_2$ are mapped to the segment $A_4K$ on the tetrahedron,
and  $Z'_1$ and $Z'_2$ are  mapped to the same point $Z$ on $A_4K$.
This  point $Z$ is  point of self intersection of the geodesic $\gamma$.

Therefore, in order that   $\gamma$  have no points of self-intersection, it is necessary that
$d > |A_4H_p|$. This implies  
\begin{equation}\label{thdcosthh}
\tanh d >  \cos(3\alpha/2)\tanh  h.
\end{equation}
The altitude $h$ of the face of the tetrahedron satisfies 
\begin{equation}\label{h}
\tanh h =\tanh a \cos\alpha/2 = \cos \alpha/2 \;\frac{\sqrt{2\cos \alpha - 1}}{\cos \alpha}.
\end{equation}
Combining  (\ref{h}) and (\ref{thdcosthh}), we obtain  
\begin{equation}\label{dist_vertex}
\tanh d>\cos\alpha/2\cos(3\alpha/2)\frac{\sqrt{2\cos\alpha-1}}{\cos\alpha},
 \end{equation}


Now we estimate the expression on the right-hand side of~(\ref{dist_vertex}) from below.
Consider  the function $\sqrt{2\cos\alpha-1}$:
$$ 2\cos\alpha-1 = 4\sin \left( \pi/6-\alpha/2 \right)  \sin \left( \pi/6+\alpha/2  \right). $$
Since the function  $\sin(\pi/6+\alpha/2)$ increases  on the interval $ (0, \pi/3)$,
then   
$$\sin(\pi/6+\alpha/2)>1/2 \;\; \textnormal{when}\;\; \alpha \in (0,  \pi/3).$$

The function $\sin \left( \pi/6-\alpha/2 \right) $ increases on the interval $ (0, \pi/3)$.
It is known that $\sin y> (2/\pi)y$ when $0<y<\pi/2$.
These imply
$$\sin(\pi/6-\alpha/2)>\frac{1}{\pi}\left( \pi/3-\alpha \right). $$
We obtain
 \begin{equation}\label{root}
\sqrt{2\cos\alpha-1} > \sqrt{ \frac{2}{3\pi} \left( \pi-3\alpha \right) }.
\end{equation}
 
The function  $\cos  (3\alpha/2)$ is decreasing for   $0<  \alpha<\pi/3$.
 It is true that $ \cos y>1- (2/\pi)y,$ when $ 0<y<\pi/2. $
Therefore,
 \begin{equation}\label{cos3alpha2}
\cos (3\alpha/2)> \frac{1}{\pi} \left( \pi-3\alpha \right).
 \end{equation}
  We  have $\cos \alpha/2 > \sqrt{3}/2$ when $0<\alpha<\pi/3$.
 
These inequalities, together with (\ref{root}) and (\ref{cos3alpha2}) give the following bound
 \begin{equation}\label{vspomd} 
\tanh d> \frac{1}{\sqrt{2\pi^3}} \left( \pi- 3\alpha \right)^{3/2}. 
 \end{equation}
The  inequality~(\ref{vspomd}) implies inequality~(\ref{distvertex2}), as required.
\end{proof}

\subsection{Uniqueness of a simple  closed geodesic  of type $(p,q)$.}
 For a regular tetrahedron in hyperbolic space
  the following analogue of Lemma~\ref{middle_spher} holds.

\begin{lemma}\label{middle_hyperbol}\textnormal{\cite{BorSuh2020}}
A simple closed geodesic on a regular tetrahedron in hyperbolic  space passes 
through the midpoints of two pairs of the opposite edges on the tetrahedron. 
\end{lemma}

\begin{proof}
Let $\gamma$ be a simple closed geodesic on a regular tetrahedron $T$
in   hyperbolic  space $\mathbb{H}^3$.
Consider  the  Cayley-Klein model of  $\mathbb{H}^3$ and
place the tetrahedron so that 
the center of circumscribed sphere of the tetrahedron  coincides with the center of the model.
Then  $T $   is represented by 
a regular tetrahedron $\tilde T$ in  Euclidean space $\mathbb{E}^3$.

A simple closed geodesic  $\gamma$  on  $T $ is represented by 
an abstract  geodesic  on $\tilde T$.
From Proposition~\ref{allgeod} we get that this  generalized   geodesic
is  equivalent to a simple closed geodesic  $\tilde\gamma$ on  $\tilde T$   in $\mathbb{E}^3$.
 From Theorem~\ref{main_th_euclid} we assume that  $\tilde\gamma$ 
passes through the midpoints of two pairs of the opposite edges on this tetrahedron. 

Label the vertices of tetrahedron $T $ and corresponding vertices of~$\tilde T$
with $A_1$, $A_2$, $A_3$ and $A_4$. 
  Suppose that $\tilde\gamma$  passes through  the midpoints  
  $\tilde X_1$, $\tilde X_2$ of the edges $A_1A_2$ and $A_3A_4$.
Consider the development of $\tilde T$ 
along  $ \tilde \gamma$ starting  from   $ \tilde X_1$.
From Corollary~\ref{sym_develop} it follows that
 this development is central-symmetric with respect to the point $\tilde  X_2$.
 
  Let $X_1$, $X_2$  be the corresponding points on 
 $\gamma$ on the edges $A_1A_2$ and $A_3A_4$ of~$T$.
 Consider the development of  $T$  onto hyperbolic plane
  along  $\gamma$ starting from the point $X_1$.
Then $\gamma$ is a line segment $X_1X'_1$  on the development.

 \begin{figure}[h]
\begin{center}
\includegraphics[width=85mm]{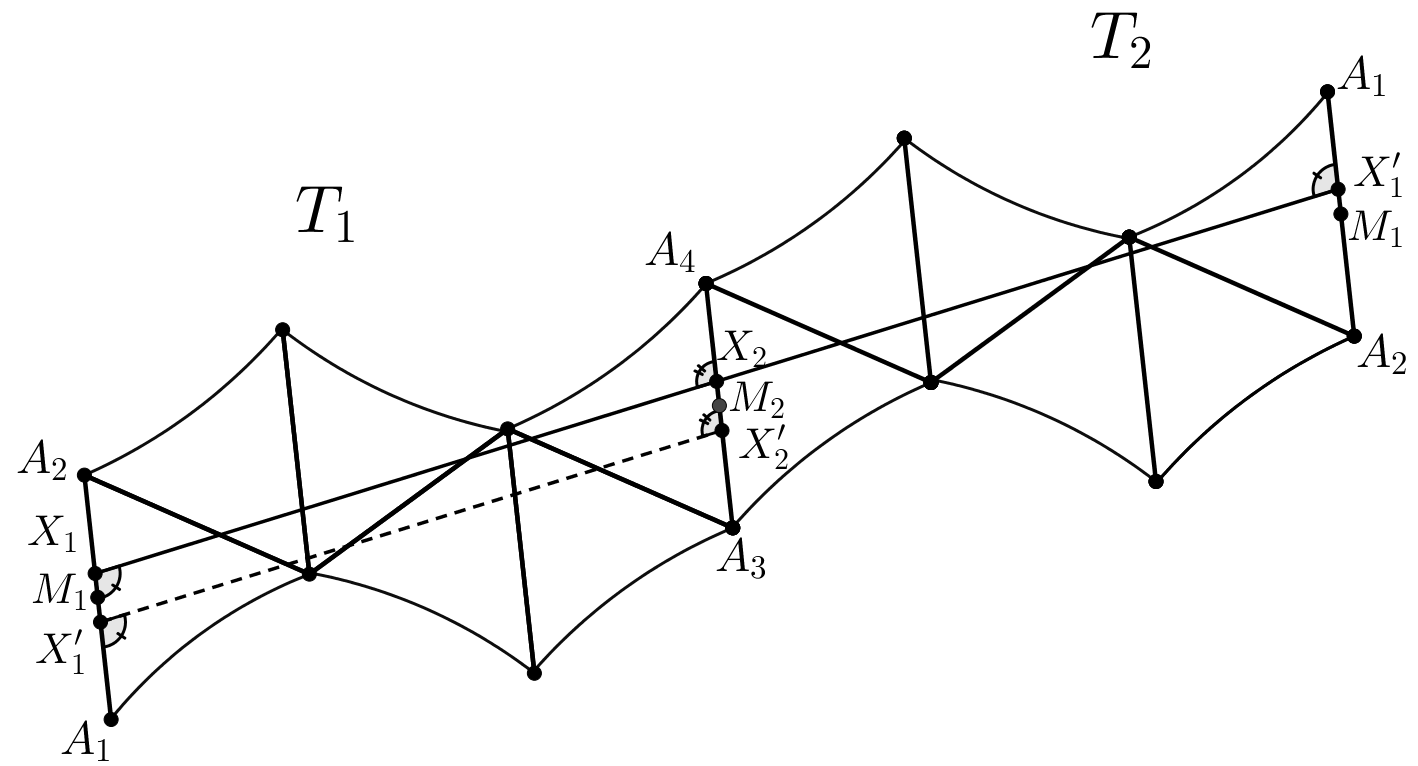}
\caption{ }
\label{newproof}
\end{center}
\end{figure}

Denote by $M_1$ and  $M_2$ the midpoints of the edges $A_1A_2$ and $A_3A_4$  respectively.
Since the rotation of  the  tetrahedron  trough  $\pi$  around  $M_1 M_2$ in hyperbolic space
 is the isometry of the  tetrahedron then,
  the development of $T$ along   $X_1X_2X'_1$ on hyperbolic plane
   is  central symmetric with the center at $M_2$.

Denote by $T_1$ and  $T_2$   the parts of the development 
along   segments $X_1X_2$ and $X_2X'_1$ respectively.
The central symmetry of the development around the point $M_2$ swap $T_1$ and $T_2$.

The edge $A_1A_2$ containing   $X'_1$ is mapped onto   $A_2A_1$ with the point  $X_1$.
Then  the point  $X'_1$  belongs to the edge $A_1A_2$ of the $T_1$, and the lengths of 
$A_2X_1$ and $X'_1A_1$ are equal.

The edge $A_3A_4$ is mapped into itself with the opposite orientation.
The point $X_2$ on  $A_3A_4$  is mapped to the point $X'_2$ on  $A_3A_4$ such that 
the  lengths of $A_4X_2$ and  $X'_2A_3$ are equal.
Moreover,  $\angle X_1X_2A_4 = \angle X'_1X'_2A_4$.
Since the geodesic is closed, then $\angle A_1X_1X_2 = \angle A_1X'_1X'_2$  (Figure \ref{newproof}).

We obtain the quadrilateral  $X_1X_2X'_2X'_1$ inside  $T_1$ the sum
of whose interior angles is $2\pi$.
Then the integral of the Gaussian curvature over
the interior of $X_1X_2X'_2X'_1$ in hyperbolic plane is zero.
This implies that the rotation takes the part $X'_2X'_1$ of the geodesic to the part $X_1X_2$. 
Hence the points $X_1$ and $X_2$ are the midpoints of the corresponding edges (Figure \ref{newproof}).

In the same way it can be proved  that  $\gamma$ passes through  the midpoints of   
 other two opposite edges on the regular tetrahedron in $\mathbb{H}^3$.
\end{proof}

\begin{corollary}\label{uniqueness_hyperbol}\textnormal{\cite{BorSuh2020}}
If two closed geodesics on the regular tetrahedron in  hyperbolic  space intersect
the edges of the tetrahedron in the same order, then they coincide. 
\end{corollary}

\subsection{Existence of a simple closed geodesic  of type $(p,q)$ on a regular tetrahedron.}

\begin{theorem}\label{existence_regtetr_hyp}\textnormal{\cite{BorSuh2020}}
On a regular tetrahedron in hyperbolic space for 
each ordered pair of coprime integers $(p, q)$, 
there exists  unique, up to the rigid motion of the tetrahedron, 
simple closed geodesic  of type $(p,q)$.
The geodesics of type $(p,q)$ exhaust all simple closed geodesics 
on a regular tetrahedron in  hyperbolic space.
\end{theorem}

\begin{proof}
 Let $\widetilde{\gamma}$  be a simple closed geodesic on 
 a regular tetrahedron $A_1A_2A_3A_4$   in Euclidean space.
 Assume, that $\widetilde{\gamma}$ passes through the midpoints $\widetilde X_1$, $\widetilde X_2$,
 $\widetilde Y_1$ and $\widetilde Y_2$ of the edges
  $A_1A_2$, $A_3A_4$, $A_1A_3$ and $A_2A_4$  respectively.
  
Consider the development  $\widetilde{T}$ of the tetrahedron along $\widetilde{\gamma}$ 
starting from the point $\widetilde X_1$ to the point $\widetilde X'_1$.
The polygon $\widetilde{T}$ consist of four equal polygons and
any two adjacent polygons can
be transformed into each other by a rotation through an angle $\pi$ around the midpoint
of their common edge - the points   $\widetilde X_2$, $\widetilde Y_1$ and $\widetilde Y_2$.
 The interior angles  of $\widetilde{T}$ are equal to
  $\pi/3$, $2\pi/3$,   $\pi$, or $4\pi/3$.
Moreover, the angle of $4\pi/3$ is obtained  if 
$\widetilde{\gamma}$ intersects three edges 
having a common vertex consecutively.

Now we take a regular triangles on  hyperbolic plane with angle  $\alpha$ at the vertices.
 Put these triangles in the same order in
which the faces of the tetrahedron were unfolded in Euclidean space along $\widetilde{\gamma}$.

 \begin{figure}[h]
\begin{center}
\includegraphics[width=75mm]{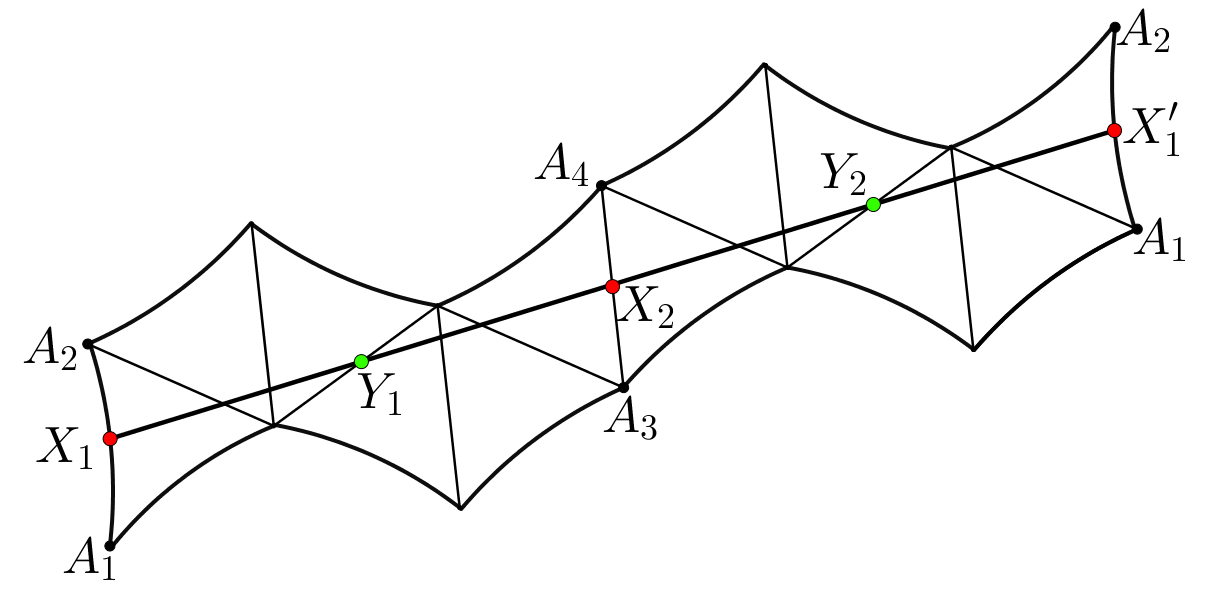}
\caption{ }
\label{dev_hyp_plane}
\end{center}
\end{figure}

In other words, we construct a polygon $T $ on hyperbolic plane that is 
formed by the same sequence of regular triangles as the polygon $\widetilde{T}$ on Euclidean plane.
Label the vertices of $T$ according to the vertices of $\widetilde{T}$.
Then the polygon $T$ corresponds to a development of a regular tetrahedron with 
the planar angle  $\alpha$ in hyperbolic space  (see Figure~\ref{dev_hyp_plane}).

Moreover, $T$ has the same property of  central symmetry   with
respect to the midpoint of the same edge as the polygon $\widetilde{T}$.
Denote by $X_1$, $X_2$, $Y_1$, $Y_2$ and $X'_1$ the midpoints 
 of the edges $A_1A_2$, $A_3A_4$, $A_1A_3$ and $A_2A_4$ of $T$  respectively.
We draw the geodesic line segment $X_1X'_1$.

By construction, the interior angles at the vertices of $T$ are equal to
$\alpha$, $2\alpha$,  $3\alpha$ or $4\alpha$,
according to the development on  Euclidean plane.

First assume that  $\alpha \in (0,  \pi/4]$.
Then  the polygon $T$ is convex
and the segment  $X_1X'_1$ lies inside $T$. 
Furthermore, $X_1X'_1$ passes through the points $X_2$, $Y_1$, $Y_2$ that are
the centers of symmetry of $T$.
 Therefore,  $X_1X'_1$ is a simple closed geodesic $\gamma$ 
on the regular tetrahedron with the planar angle  $\alpha \in (0, \pi/4]$
 in hyperbolic space.

Now we increase the angle $\alpha$ starting from $\alpha =\pi/4$.
Then the polygon $T$ is not convex because it contains the interior angles $4\alpha>\pi$.

Let $\alpha_0$ be the supremum of  $\alpha$ 
for which the segment  $ X_1X_2$ lies inside $T$.
Suppose $\alpha_0<\pi/3$.
For all $\alpha<\alpha_0$ the segment $X_1X'_1$ 
lies entirely inside $T$ and 
 it is a simple closed geodesic  $\gamma$  on 
the regular tetrahedron in $\mathbb{H}^3$.
The distance $d$  from  the  vertices of the  tetrahedron to $\gamma$
satisfies~(\ref{distvertex2}).
Therefore there exists 
$\alpha_1=\alpha_0+\varepsilon$   such that  
 the segment $ X_1X_2$ lies entirely inside $T$.
This contradicts the maximality of $\alpha_0 $.
Thus  $\alpha_0=\pi/3$.

It follows that for any $\alpha \in (0, \pi/3)$ there is a simple closed geodesic of type $(p,q)$
 on a regular  tetrahedron with a planar angle $\alpha$ in hyperbolic space.

From Corollary~\ref{uniqueness_hyperbol} it follows 
the uniqueness of a simple closed geodesic of type $(p,q)$
on a regular tetrahedron in $\mathbb{H}^3$.
This geodesic has $p$ points on  each of two opposite edges of the tetrahedron,
$q$ points on  each of another two opposite edges,
and $(p+q)$ points on each edge of the third pair of opposite edges.
For any coprime integers $(p, q)$, $0 \le p<q$,
 there exist three simple closed geodesic of type $(p,q)$
on a regular tetrahedron in $\mathbb{H}^3$.
They coincide by the rotation of the tetrahedron
by the angle $2\pi/3$ or $4\pi/3$ about the altitude constructed from a vertex to the opposite face.

Since any simple closed geodesic on a regular tetrahedron in $\mathbb{H}^3$ is equivalent to 
a simple closed geodesic on a regular tetrahedron in $\mathbb{E}^3$, then there are not  other 
simple closed geodesics on a regular tetrahedron in $\mathbb{H}^3$.
\end{proof}

\subsection{Existence of a simple closed geodesic   of type $(p,q)$ on a generic tetrahedron.}
In Euclidean space $\mathbb{E}^3$, there is no simple closed geodesic  on a generic tetrahedron.
Protasov~\cite{Pro07} gave an upper bound for the number of simple closed geodesics
depending on the largest deviation from $\pi$ of the sum of planar angles at the vertices of
the tetrahedron .
The situation in hyperbolic space is quite different provided that the planar
angles of the tetrahedron are sufficiently small. 
Borisenko proved the following result.

\begin{theorem}\label{gener_tetr_hyper}\textnormal{\cite{Bor2022}}
If the planar angles of a tetrahedron in hyperbolic space are at most $\pi/4$,
then for any pair of coprime natural numbers $(p, q)$ there exist  
three simple closed geodesics of type $(p, q)$,
disregarding isometries of the tetrahedron.
\end{theorem}

\begin{proof}
 Let $\widetilde{\gamma}$  be a simple closed geodesic on 
 a regular tetrahedron $A_1A_2A_3A_4$ in Euclidean space.
Consider the development  $\widetilde{T}$ of the tetrahedron along $\widetilde{\gamma}$ 
starting from the point $\widetilde X_1$ on $A_1A_2$ to the point $\widetilde X'_1$.

Consider a generic tetrahderon in hyperbolic space. 
For convenience we can also label the tetrahedron's vertices with $A_1$, $A_2$, $A_3$ and $A_4$.
Develop this tetrahedron onto the hyperbolic plane in the same order as 
the development $\widetilde{T}$ is unfolded, starting from the edge $A_1A_2$.

As it was shown in  the proof of Theorem~\ref{existence_regtetr_hyp},
at most four facets can meet at one vertex of the development.
Hence if $\alpha \le \pi/4$, then the development is a convex polygon. 

However, there are at most two facets meeting at each of the vertices $A_1$, $A_2$, $A'_1$, $A'_2$,
where  $A_1A_2$ is starting  edge and $A'_1A'_2$ is finishing.
Therefore the angles at these vertices are at most $\pi/2$.

 Consider the quadrilateral $A_1A_2A'_2A'_1$.
 Take points $X(s)$ on $A_1A_2$ and $X'(s)$ on $A'_1A'_2$
such that $X(0) = A_1$, $X'(0) = A'_1$, and the lengths of $A_1X(s)$ and $A'_1X'(s)$ are both
equal to $s$ (Figure~\ref{hyperbol_case_general}).

 \begin{figure}[h]
\begin{center}
\includegraphics[width=85mm]{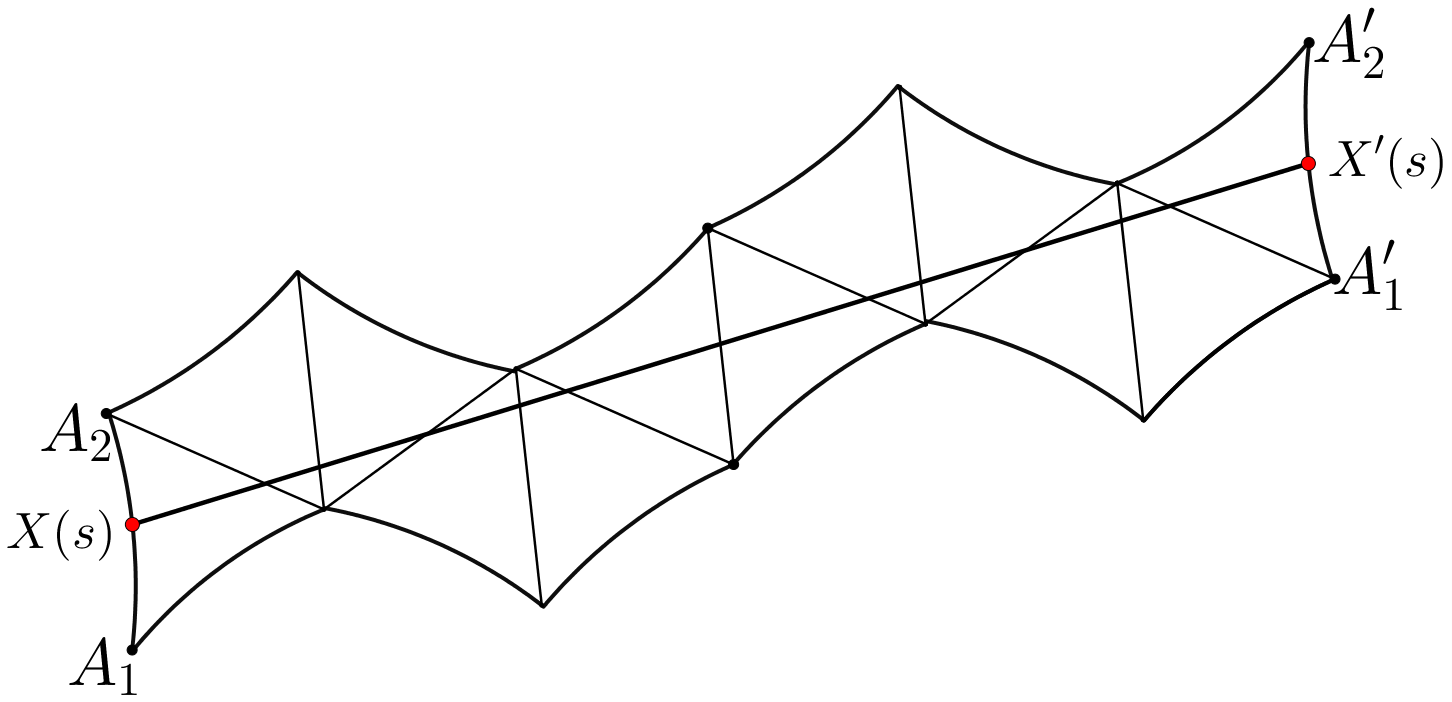}
\caption{ }
\label{hyperbol_case_general}
\end{center}
\end{figure}

 For $s = 0$, the sum of the angles $\angle A_1$ and $\angle A'_1$ measured from inside the polygon
is less than $\pi$. For $s = |A_1A_2|$, 
the sum of  $\angle A_2$ and  $\angle A'_2$ measured from outside the 
polygon is greater than  $\pi$.
 Therefore, there is $s_0$ such that the sum of  $\angle X(s_0)$ and $\angle X'(s_0)$ equals $\pi$.
 The line segment $X(s_0)X'(s_0)$ on the development corresponds
to a simple closed geodesic of type $(p, q)$ on the tetrahedron in $\mathbb{H}^3$.

Since for any ordered pair of coprime integers $(p, q)$ 
 there exist three simple closed geodesic of type $(p,q)$
on a regular tetrahedron in $\mathbb{E}^3$, disregarding isometries of the tetrahedron, 
then similarly we can construct three simple closed geodesic of type $(p,q)$
on a tetrahedron in $\mathbb{H}^3$ with planar angle at most $\pi/4$.
\end{proof}

\subsection{The number of simple closed geodesics.}
Let $N(L, \alpha) $ be a number of simple closed geodesics of length  not greater than $L$ 
on a regular tetrahedron  with  planar angle $\alpha$ in hyperbolic space.
In~\cite{BorSuh2020} it was shown that 
$$ N(L, \alpha) = c(\alpha) L^2 +O(L\ln L),$$
where $O(L\ln L) \le CL\ln L$ when $L \rightarrow +\infty$, and
\begin{equation}\label{constant_old}
 c(\alpha) = \frac{27}{ 16 
 \left( 
  \ln 
  \frac{1 - \frac{\sqrt{3}}{2} 
 \left( 1 - \frac{3\alpha}{\pi} \right)^3 \left(1- \frac{\alpha^2}{4}\right)  }
{1 - \frac{\sqrt{3}}{2} \left( 1 - \frac{3\alpha}{\pi} \right)^3 \left(1+ \frac{\alpha^2}{4}\right) }
+
\ln   \frac{1 + \frac{\sqrt{3}}{4}\left( 1 - \frac{3\alpha}{\pi}  \right)}
{1 - \frac{\sqrt{3}}{4}\left( 1 - \frac{3\alpha}{\pi} \right)}   
\right) ^2}, \notag
\end{equation}
 \begin{equation}\label{cpi3}
\lim_{\alpha \rightarrow \frac{\pi}{3}} c(\alpha) = +\infty;
\; \; \; \; 
\lim_{\alpha \rightarrow 0} c(\alpha) =
 \frac{27}{ 16\left(
   \ln  
    \frac{1 + \frac{\sqrt{3}}{4} }
{1 - \frac{\sqrt{3}}{4} }    \right) ^2}.
\end{equation}
 
This result was proved using Proposition~\ref{gengeodstr} about the structure of a simple closed geodesic 
on a regular tetrahedron. 
 
 In current paper we   improve the constant $c(\alpha)$ 
using the estimations obtained in~\cite{BorSuh2020}.

\begin{lemma}\label{condpq}
If the length of a  simple closed geodesic  of type $(p,q)$ 
on a regular tetrahedron in hyperbolic space is not greater than $L$, then 
  \begin{equation} 
L  \ge 2(p+q) \ln \left( 2\sqrt{3} \left( 1-\frac{3\alpha}{\pi}\right) +1  \right). \notag
\end{equation}
 where $\alpha$ is the plane angle  of a face of the tetrahedron.
\end{lemma}
 
\begin{proof}
Let $\gamma$ be a simple closed geodesic of type $(p,q)$, $0\le p <q$,
on a regular tetrahedron $A_1A_2A_3A_4$
in hyperbolic space. 

Assume that $\gamma$ has $q$ points on the edges $A_1A_2$ and $A_3A_4$,
$p$ points on $A_1A_4$ and $A_2A_3$ and $p+q$ points on $A_2A_4$ and $A_1A_3$. 
Denote by $B_1, \dots, B_{p+q}$ points of  $\gamma$  on $A_1A_3$  and
by $B'_1, \dots, B'_{p+q}$ points of  $\gamma$  on $A_2A_4$.

Consider the development of the faces $A_3A_1A_4$ and $A_1A_4A_2$ onto the plane. 
The geodesic segment starting at the point $B_i$, where $i=1, \dots, p$,
 goes through the edge $A_1A_4$ to the point $B'_{q+i}$.
 Similarly on the development of the faces $A_1A_2A_3$ and $A_2A_3A_4$ 
there are $p$ segments of $\gamma$ connecting $B'_i$ and $B_{q+i}$, $i~=~1,~\dots,~p$ and 
passing through the edge $A_2A_3$. 
 
 On the faces $A_4A_1A_2$ and $A_1A_2A_3$ the geodesic segments
  $B_i B'_{q-(i-1)}$, $i~=~1,~\dots,~q$, pass through the  edge $A_1A_2$. 
 Similarly on the development of the faces $A_2A_4A_3$ and $A_4A_3A_1$ 
 there are $q$ geodesic segments $B_{p+i}B'_{(p+q)-(i-1)}$, $i~=~1,~\dots,~q$ 
 (see Figure~\ref{structure_of_geod}).

 \begin{figure}[h]
\begin{center}
\includegraphics[width=75mm]{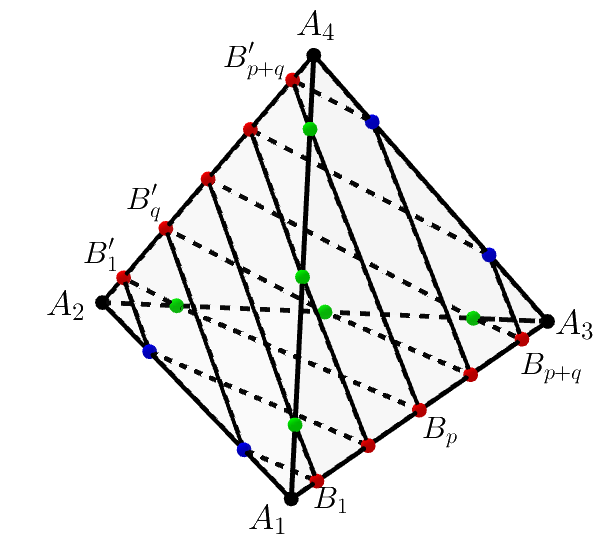}
\caption{ }
\label{structure_of_geod}
\end{center}
\end{figure}

 Therefore the geodesic $\gamma$ consists of $2(p+q)$ segments, that 
 connect opposite edges of the tetrahedron. 
Let us evaluate from below the length of these segments. 
Consider the quadrilateral obtaining by unfolding of the faces
$A_2A_1A_4$ and $A_1A_4A_3$.
The minimum distance between points on the edges $A_2A_4$  and $A_1A_3$
  is achieved at $H_1H_2$ perpendicular to these edges. 
  Since the planar angle of the tetrahedron $\alpha < \pi/3$, then $H_1H_2$ 
  lies inside the  quadrilateral $A_3A_1A_4A_2$ and passes through the 
  midpoint $M$ of the edge $A_1A_4$ (see Figure~\ref{hyper_dev_two_edg}). 
 \begin{figure}[h]
\begin{center}
\includegraphics[width=75mm]{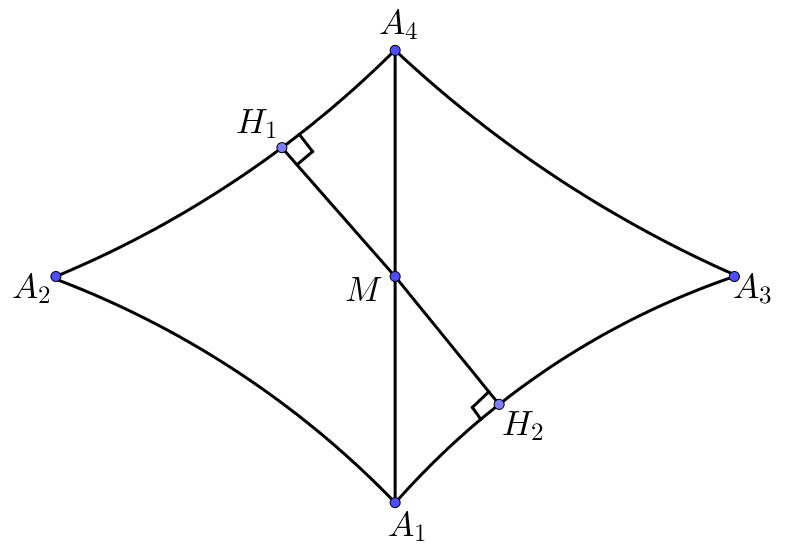}
\caption{ }
\label{hyper_dev_two_edg}
\end{center}
\end{figure}
From the triangle $A_4MH_1$  we have
$$ \sinh |MH_1| =  \sinh (a/2) \sin \alpha $$
Using (\ref{a_hyp}) we get 
\begin{equation} 
 \sinh |MH_1| = \cos(\alpha/2)  \sqrt{2\cos \alpha-1}. \notag
\end{equation}
Using $$ 2\cos \alpha-1 = \frac{\cos(3\alpha/2)}{\cos(\alpha/2)} $$ we get
\begin{equation} 
 \sinh |MH_1| =   \sqrt{ \cos(\alpha/2)\cos(3\alpha/2) }. \notag
\end{equation}
The inequality~(\ref{cos3alpha2}) together with  $\cos \alpha/2 > \sqrt{3}/2$ implies
\begin{equation} 
 \sinh |MH_1| \ge   \sqrt{ \frac{\sqrt{3}}{2}  \left( 1-\frac{3\alpha}{\pi}\right)  }.  
\end{equation}
Consider the function  arsinh($x$):
\begin{equation} 
2\textnormal{arsinh(x)} = 2 \ln \left( x + \sqrt{x^2+1} \right) = \ln (2x^2+1+2x\sqrt{x^2+1})
> \ln (4x^2+1). \notag
 \end{equation}
 This implies
$$   |H_1H_2| \ge \ln \left(2\sqrt{3}  \left( 1- 3\alpha/\pi\right) +1\right). $$
We obtain that the length $L $ of a  simple closed geodesic $\gamma$ of type $(p,q)$ satisfies
  \begin{equation} 
L  \ge 2(p+q)  \ln \left(2\sqrt{3} \left( 1- 3\alpha/\pi\right) +1\right). \notag
\end{equation}
\end{proof}

Euler's  function $\phi(n)$ is equal to the number of
positive integers not greater than  $n $ and prime to $n \in \mathbb N$.
From \cite[Th. 330.]{Hard} we know
\begin{equation}\label{phi}
 \sum\limits_{n=1}^{x}  \phi(n) = \frac{3}{\pi^2}x^2+O(x\ln x),
\end{equation}
where $O(x\ln x) < C x\ln x $, when $x \rightarrow +\infty$.

 Denote by $\psi(x)$ the number of  pairs of coprime integers $(p,q)$
  such that $p<q$ and $p+q\le x,$ $x\in \mathbb R$.
  Suppose  $\hat \psi(y)$  is equal to the number of  pairs
   of coprime integers $(p,q)$ such that $p<q$ and 
$p+q= y$, $y\in \mathbb N$.
From the definitions we get
\begin{equation}\label{psihatpsi}
\psi(x) = \sum\limits_{y=1}^{x} \hat \psi(y) 
\end{equation}

If $(p, q)=1$ and $p+q=y$, then $(p, y)=1$ and $(q, y)=1$.
 Consider Euler's  function $\phi(y)$.
We obtain  that the set of  integers not greater than and prime to $y$ are separated
 into the  pairs of coprime integers  
$(p,q)$ such that $p<q$ and $p+q= y$
 It follows that  $\phi(y)$ is even  and $\hat \psi(y) = \phi(y)/2$.
From (\ref{psihatpsi}) we have
\begin{equation}
\psi(x) = \frac{1}{2} \sum\limits_{y=1}^{x} \phi(y).\notag
\end{equation}

The formula (\ref{phi}) implies
 \begin{equation}\label{psihatpsi2}
\psi(x) = \frac{3}{2\pi^2}x^2+O(x\ln x),
\end{equation}
where $O(x\ln x) < C x\ln x $ when $x \rightarrow +\infty$.
 
 Using this asymptotic it can be proved following result.

 \begin{theorem}\label{estimate_number}
Let $N(L, \alpha) $ be the number of simple closed geodesics of length  not greater than $L$ 
on a regular tetrahedron  with  plane angles of the faces equal to $\alpha$ in hyperbolic space.
Then  
\begin{equation}\label{NLalpha}
 N(L, \alpha) = c(\alpha) L^2 +O(L\ln L),
\end{equation}
where
\begin{equation}
 c(\alpha) = \frac{9}{8\pi^2 
 \ln \left(2\sqrt{3} \left( 1- 3\alpha/\pi\right) +1\right)   },  \notag
\end{equation}
 \begin{equation}\label{cpi3}
\lim_{\alpha \rightarrow \frac{\pi}{3}} c(\alpha) = +\infty;
\; \; \; \; 
\lim_{\alpha \rightarrow 0} c(\alpha) =
 \frac{9}{ 8 \pi^2 \ln \left(2\sqrt{3}  +1\right)}.
\end{equation}
and  $O(L\ln L) \le CL\ln L$ when $L \rightarrow +\infty$,
 \end{theorem}
 
\begin{proof}
To each ordered pair of coprime integers $(p,q)$, $p<q$  there correspond 
three different geodesics on the regular tetrahedron.
 We have
\begin{equation}
N(L, \alpha) =3\psi \left(   
\frac { L }{2\ln \left(2\sqrt{3} \left( 1- 3\alpha/\pi \right) +1\right) }   \right) \notag
\end{equation}
Using  (\ref{psihatpsi2}),  we get 
\begin{equation}
N(L, \alpha) =\frac{9}{8\pi^2 
 \ln \left(2\sqrt{3} \left( 1- 3\alpha/\pi\right) +1\right)   }
L^2   +O(L\ln L), \notag
\end{equation}
 when $L \rightarrow +\infty$. 
\end{proof}

  In  work~\cite{Rivin} of Rivin it was shown that 
  for any hyperbolic structure  on a sphere  with $n$ boundary components,
  the number of simple closed geodesics of length bounded by $L$  on it
     grows like $L^{2n-6}$ as $L~\rightarrow~\infty$.

  From Lemma~\ref{nesdist} we know that  
the  distances from the vertices of the tetrahedron to 
a simple closed geodesic   is greater then  $d(\alpha)$,
 where $d_0(\alpha)$ is from~(\ref{distvertex2}).
  This estimation  holds also for a generic tetrahedron in hyperbolic space. 

We can consider tetrahedron as a non-compact surface with  
 regular Riemannian metric of constant negative curvature with $4$ boundary components.
 From Lemma~\ref{intersect_four_edges} it follows that there is no   simple closed geodesic 
 that is  boundary parallel.
 From~(\ref{NLalpha}) we get, that  the number  $N(L, \alpha)$ of  
 simple closed geodesics of length $\le L$ on a regular  tetrahedron
  is asymptotic to $L^2$ when $L \rightarrow +\infty$.

If the planar angle $\alpha$  of the tetrahedron  goes to zero, 
then the vertices of the tetrahedron tend  to infinity.
The limiting surface is  homeomorphic to a sphere with $4$ cusps .
It is shown in~\cite{McShane_Rivin} that any cusp on hyperbolic surface  has a neighborhood 
bounded by  horocycle  curve of length $2$. 
There is no simple closed geodesic intersecting this neighborhood. 
 Then the asymptotic of the number of 
  simple closed geodesics  on a sphere  with $n$ cusps 
 is equal to the  asymptotic of the  number 
 of simple closed geodesics on a sphere  with $n$ boundary components.


B.Verkin Institute for Low Temperature Physics and Engineering 
of the National Academy of Sciences of Ukraine,
Kharkiv, 61103, Ukraine

Mathematisches Institut, WWU M\"unster, Einsteinstrasse 62, D-48149, M\"unster

\textsc {\textit{E-mail address:}}  suhdaria0109@gmail.com, s.darya@uni-muenster.de

\end{document}